\documentclass[11pt]{article}
\usepackage{amsfonts, amsbsy, amssymb, eucal}
\usepackage{graphicx}

\begin{document}

%%%%%%%%%%%%%%%%%%%

%LETTERS
\newcommand{\eps}{\varepsilon}
\newcommand{\seps}{\sqrt{\eps}}
\def\omeps{{\omega \over \seps}}
\def\omepsi{{\omega \over \eps}}

%bold
\newcommand{\x}{\boldsymbol x}
\newcommand{\s}{\boldsymbol s}
\newcommand{\io}{\boldsymbol \iota}
\newcommand{\p}{\boldsymbol p}
\newcommand{\q}{\boldsymbol q}
\newcommand{\g}{\boldsymbol g}
\newcommand{\ar}{\boldsymbol a}
\newcommand{\br}{\boldsymbol b}
\newcommand{\ur}{\boldsymbol u}
\newcommand{\ka}{\boldsymbol k}
\newcommand{\w}{\boldsymbol \omega}
\newcommand{\ps}{\boldsymbol \psi}
\newcommand{\xii}{\boldsymbol \xi}

%frank
\newcommand{\fb}{\mathfrak b}
\newcommand{\fB}{\mathfrak B}
\newcommand{\fC}{\mathfrak C}
\newcommand{\fD}{\mathfrak D}
\newcommand{\fN}{\mathfrak N}
\newcommand{\fW}{\mathfrak W}
\newcommand{\fS}{\mathfrak S}
\newcommand{\fI}{\mathfrak I}
\newcommand{\fri}{\mathfrak i}
\newcommand{\fp}{\mathfrak{p}}

%small
\newcommand{\de}{{\scriptstyle \Delta}}
\newcommand{\te}{{\scriptstyle T}}
\newcommand{\ie}{{\scriptstyle I}}

%misc
\newcommand{\ld}{,\ldots,}
\newcommand{\spr}[2]{\langle#1,#2\rangle}
\newcommand{\avr}[1]{\langle#1\rangle}
\newcommand{\ti}[1]{{}^{\mathfrak t}{#1}^{\mathtt{-1}}}
\def\={\stackrel{\rm def}{=}}

%standard
\def\L{{\mathbb L}}
\def\B{{\mathbb B}}
\def\R{{\mathbb R}}
\def\Rn{{\mathbb R}^n}
\def\Rnp{{\mathbb R}^{n+1}}
\def\T{{\mathbb T}}
\def\N{{\mathbb N}}
\def\M{{\mathbb M}}
\def\C{{\mathbb C}}
\def\Z{{\mathbb Z}}
\def\Q{{\mathbb Q}}

%operators
\def\Dlw{D_{\lambda_0,\omega}}
\def\Dlwl{D_{\lambda_0,\omega,\Lambda}}
\newcommand{\Diff}{\lambda D_s+\spr{\omega }{D_\varphi} }
\newcommand{\Difff}{\lambda D_s+\spr{\omega }{D_\varphi} + \spr{\Omega z
}{D_z}}

%manifolds
\newcommand{\It}{{\cal I}\times\T^n}
\newcommand{\cyl}{\cal C}
\newcommand{\cyle}{\mathfrak C}
%sets
\newcommand{\di}{{\fW}^n_{\tau,\gamma}}
\newcommand{\defp}{{\fD}_{\fp,\de}({\cal C})}
\newcommand{\defpit}{{\fD}_{\fp,\de}(\It)}
\newcommand{\defppi}{{\fD}_{\fp,\de}(\Pi\times\T^n)}

%spaces
\newcommand{\bst}{{\fB}_\sigma(\T^n)}
\newcommand{\bset}{{\fB}_{\kappa,\sigma}(\B^m\times\T^n)}
\newcommand{\bstprime}{{\fB}_{\sigma'}(\T^n)}
\newcommand{\bsetprime}{{\fB}_{\kappa,\sigma'}(\B^m\times\T^n)}
\newcommand{\bstk}{{\fB}_{\kappa,\sigma}(T^*\T^n)}

\newcommand{\bsitb}{{\fB}_{\beta,\fp}({\cal I}\times\T^n)}
\newcommand{\bsitbj}{{\fB}_{\beta,\fp}^j({\cal I}\times\T^n)}
\newcommand{\bsit}{{\fB}_{\fp}({\cal I}\times\T^n)}
\newcommand{\bsitj}{{\fB}_{\fp}^j({\cal I}\times\T^n)}
\newcommand{\bsitprime}{{\fB}_{\fp'}({\cal I}\times\T^n)}

\newcommand{\bskit}{{\fB}_{\kappa,\mathfrak{p}}[T^*(\It)]}
\newcommand{\bskitprime}{{\fB}_{\kappa',\mathfrak{p}'}[T*(\It)]}

\newcommand{\bs}{{\fB}_{\mathfrak{p}}({\cal C})}
\newcommand{\bsone}{{\fB}^1_{\mathfrak{p}}({\cal C})}
\newcommand{\bstwo}{{\fB}^2_{\mathfrak{p}}({\cal C})}
\newcommand{\bsj}{{\fB}^j_{\mathfrak{p}}({\cal C})}

\newcommand{\bsprime}{{\fB}_{\mathfrak{p}'}({\cal C})}
\newcommand{\bsprimeb}{{\fB}_{\mathfrak{p}'}({\cal C_\beta})}
\newcommand{\bsoneprime}{{\fB}^1_{\mathfrak{p}'}({\cal C})}
\newcommand{\bstwoprime}{{\fB}^2_{\mathfrak{p}'}({\cal C})}
\newcommand{\bsjprime}{{\fB}^j_{\mathfrak{p}'}({\cal C})}

\newcommand{\bsm}{{\fB}^{\wedge}_{\mathfrak{p}}({\cal C})}
\newcommand{\bsmprime}{{\fB}_{\mathfrak{p}'}^{\wedge}({\cal C})}

\newcommand{\bsg}{{\fB}^{(n,\wedge)}_{\mathfrak{p}}({\cal C})}
\newcommand{\bsgprime}{{\fB}^{(n,\wedge)}_{\mathfrak{p}'}({\cal C})}

\newcommand{\bsz}{{\fB}^{(0,n+1)}_{\fp}({\cal C})}
\newcommand{\bszprime}{{\fB}^{(0,n+1)}_{\fp'}({\cal C})}

\newcommand{\bsk}{{\fB}_{\kappa,\mathfrak{p}}(T^*{\cal C})}
\newcommand{\bskz}{{\fB}_{\kappa,\mathfrak{p}}(T^*{\cal C}\times
\B^{2(m-1)})}

\newcommand{\bsbk}{{\fB}_{\beta,\kappa,\mathfrak{p}}(T^*{\cal C})}
\newcommand{\bsbkz}{{\fB}_{\beta,\kappa,\mathfrak{p}}(T^*{\cal C}\times
\B^{2(m-1)})}
\newcommand{\bskprime}{{\fB}_{\kappa',\mathfrak{p}'}(T^*{\cal C})}
\newcommand{\bskprimez}{{\fB}_{\kappa',\mathfrak{p}'}(T^*{\cal C}
\times \B^{2(m-1)})}

\newcommand{\bsbkprime}{{\fB}_{\beta,\kappa',\mathfrak{p}'}(T^*{\cal C})}
\newcommand{\bsbkprimez}{{\fB}_{\beta,\kappa',\mathfrak{p}'}(T^*{\cal C}
\times \B^{2(m-1)})}

\newcommand{\bscheck}{ {\fB}_{\fp} (\check{\cal C}) }
\newcommand{\bsinf}{{\fB}_{r,\infty,\rho,\sigma}({\cal C})}
\newcommand{\bskinf}{{\fB}_{r,\infty,\rho,\sigma}({\cal C})}
\newcommand{\bskcheck}{ {\fB}_{\kappa;\fp} (T^*\check{\cal C}) }
\newcommand{\bshat}{ {\fB}_{\fp} (\hat{\cal C}) }
\newcommand{\bshatprime}{ {\fB}_{\fp'} (\hat{\cal C}) }
\newcommand{\bskhat}{ {\fB}_{\kappa;\fp} (T^*\hat{\cal C}) }

\newcommand{\bspi}{\fB_{\fp}(\Pi\times\T^n)}
\newcommand{\bspiprime}{\fB_{\fp'}(\Pi\times\T^n)}

\newcommand{\bsp}{{\fB}_\sigma(\T^n)\times[\bs]^{n+1}}
\newcommand{\bspprime}{{\fB}_{\sigma'}(\T^n)\times[\bsprime]^{n+1}}

%more stuff for whiskers
\newcommand{\bse}{{\fB}_{\mathfrak{p}}({\mathfrak C})}
\newcommand{\bseone}{{\fB}^1_{\mathfrak{p}}({\mathfrak C})}
\newcommand{\bsetwo}{{\fB}^2_{\mathfrak{p}}({\mathfrak C})}
\newcommand{\bsej}{{\fB}^j_{\mathfrak{p}}({\mathfrak C})}

\newcommand{\bseprime}{{\fB}_{\mathfrak{p}'}({\mathfrak C})}
\newcommand{\bseprimeb}{{\fB}_{\mathfrak{p}'}({\cal C_\beta})}
\newcommand{\bseoneprime}{{\fB}^1_{\mathfrak{p}'}({\mathfrak C})}
\newcommand{\bsetwoprime}{{\fB}^2_{\mathfrak{p}'}({\mathfrak C})}
\newcommand{\bsejprime}{{\fB}^j_{\mathfrak{p}'}({\mathfrak C})}

\newcommand{\bsem}{{\fB}^{\wedge}_{\mathfrak{p}}({\mathfrak C})}
\newcommand{\bsemprime}{{\fB}_{\mathfrak{p}'}^{\wedge}({\mathfrak C})}

\newcommand{\bseg}{{\fB}^{(n,\wedge,m)}_{\mathfrak{p}}({\mathfrak C})}
\newcommand{\bsegprime}{{\fB}^{(n,\wedge,m)}_{\mathfrak{p}'}({\mathfrak C})}

\newcommand{\bsek}{{\fB}_{\kappa,\mathfrak{p}}(T^*{\mathfrak C})}

\newcommand{\bsebk}{{\fB}_{\beta,\kappa;\mathfrak{p}}(T^*{\mathfrak C})}

\newcommand{\bsekprime}{{\fB}_{\kappa';\mathfrak{p}'}(T^*{\mathfrak C})}

\newcommand{\bsebkprime}{{\fB}_{\beta,\kappa';\mathfrak{p}'}(T^*{\mathfrak C})}

\newcommand{\bsecheck}{ {\fB}_{\fp} (\check{\mathfrak C}) }
\newcommand{\bseinf}{{\fB}_{r,\infty,\rho,\sigma}({\mathfrak C})}
\newcommand{\bsekinf}{{\fB}_{r,\infty,\rho,\sigma}({\mathfrak C})}
\newcommand{\bsekcheck}{ {\fB}_{\kappa;\fp} (T^*\check{\mathfrak C}) }
\newcommand{\bsehat}{ {\fB}_{\fp} (\hat{\mathfrak C}) }
\newcommand{\bsehatprime}{ {\fB}_{\fp'} (\hat{\mathfrak C}) }
\newcommand{\bsekhat}{ {\fB}_{\kappa;\fp} (T^*\hat{\mathfrak C}) }

\newtheorem{definition}{Definition}
\newtheorem{assumption}{Assumption}
\newtheorem{hyp}{Hypothesis}
\newtheorem{theorem}{Theorem}
\newtheorem{lemma}{Lemma}[section]
\newtheorem{corollary}{Corollary}[theorem]
\newtheorem{proposition}{Proposition}

\newcounter{pups}[section]
\setcounter{pups}{0}

\newenvironment{remark}{

\medskip\noindent\addtocounter{pups}{1}{\tt Remark \arabic{section}.\arabic{pups}:}}{

\medskip\noindent}

\renewcommand{\theenumi}{\roman{enumi}}

%%%%%%%%%%%%%%%%%%%
\title{A model for separatrix splitting near multiple resonances
\footnotetext{AMS subject classification 70H08, 70H20}}
\author{M. Rudnev\thanks{Contact address: Department of Mathematics,
University of Bristol, University Walk, Bristol BS8 1TW, UK;\newline
e-mail: {\tt m.rudnev@bris.ac.uk}} $\,$and V. Ten\thanks{Contact
address: Department of Mathematics, University of Bristol,
University Walk, Bristol BS8 1TW, UK;\newline  e-mail: {\tt
v.ten@bris.ac.uk} }}
\date{January 13, 2004}
\maketitle

\begin{abstract}
We propose a model for local dynamics of a perturbed convex
real-analytic Liouville-integrable Hamiltonian system near a
resonance of multiplicity $1+m,\,m\geq 0$. Physically, the model
represents a toroidal pendulum, coupled with a Liouville-integrable
system of $n$ non-linear rotators via a small analytic potential.
The global bifurcation problem is set-up for the $n$-dimensional
isotropic manifold, corresponding to a specific homoclinic orbit of
the toroidal pendulum. The splitting of this manifold can be
described by a scalar function on an $n$-torus, whose $k$th Fourier
coefficient satisfies the estimate
$$O\left(e^{-\,\rho|k\cdot\omega| -
|k|\sigma}\right),\;\,k\in\Z^n\setminus\{0\},$$ where
$\omega\in\R^n$ is a Diophantine rotation vector of the system of
rotators; $\rho\in(0,{\pi\over2})$ and $\sigma>0$ are the
analyticity parameters built into the model. The estimate, under
suitable assumptions would generalize to a general multiple
resonance normal form of a convex analytic Liouville integrable
Hamiltonian system, perturbed by $O(\eps)$, in which case
$\omega_j\sim\omeps,\,j=1,\ldots,n.$
\end{abstract}

\subsection*{1. Introduction and main
result}\addtocounter{section}{1}
\renewcommand{\theequation}{\arabic{section}.\arabic{equation}}

The main objective of this paper is to create a template to extend
the theory for exponentially small separatrix splitting in Liouville
near-integrable Hamiltonian systems near simple resonances, i.e.
resonances of multiplicity $1$, to the case of multiple resonances,
of multiplicity $1+m$, $m\geq0$. The interest in such a theory is
dictated by the fact that the normal form theory and Nekhoroshev
theorem, resulting in exponentially long time stability (see e.g.
\cite{LMS}) is well  developed for resonances of all multiplicities,
whereas the exponentially small splitting phenomenon, resulting in
similar exponents, has been quantitatively studied so far only in
the special case of multiplicity one resonances.

It is well known (see e.g. \cite{DS}) that a convex analytic
Liouville near-integrable Hamiltonian system, with Hamiltonian
\begin{equation}
H(\p,\q)=H_0(\p)+\varepsilon H_1(\p,\q), \label{nli}\end{equation}
where $(\p,\q)\in\R^{n+1+m}\times \T^{n+1+m}$ are the action-angle
variables on $T^*\T^{n+1+m}$, where $\T=\R/2\pi\Z$, can be localized
in the action space near a multiplicity $1+m$ resonant action value
$\p_0$. Namely suppose $\p_0$ is such that the kernel of the scalar
product $\spr{DH_0(\p_0)}{\ka},\,\ka\in\Z^{n+1+m}$ is some $1+m$
dimensional sublattice in $\Z^{n+1+m}$. In this case, without loss
of generality one can render $DH_0(\p_0)=(\omega,0)\in\R^{n+1+m}.$
In addition, we further assume that $\omega\in\R^n$ is Diophantine,
i.e.
\begin{equation}
\forall\,k\in\Z^{n}\setminus\{0\},\;\;
|\spr{k}{\omega}|\geq\vartheta|k|^{-\tau_n},\label{dio}
\end{equation}
for some $\vartheta>0$ and $\tau_n\geq n-1$ (for $n=1$ this
obviously boils down to $\omega\neq0$).

After a canonical change of variables, preserving the phase space
bundle structure and time scaling the Hamiltonian (\ref{nli}) can be
cast into the following normal form:

\begin{equation}
H_{\rm nf}(\p,\q)\;=\;\spr{{\omega\over\seps}}{\iota}+
\,{1\over2}\spr{\p}{Q_{\rm nf}\p} + U(x_0,\ldots,x_m)\,\,+\;\; [
f_{\rm nf}(\q) + \spr{\p}{\g_{\rm nf}(\p,\q)}], \label{nfh}
\end{equation} where $\p=(\iota,y_0,y_1,\ldots,y_m)\in\R^{n+1+m},\,\q=(\varphi,x_0,x_1,\ldots
x_m)\in\T^{n+1+m}$, $Q_{\rm nf}$ is a constant symmetric matrix, and
the pair $(f_{\rm nf},\g_{\rm nf})=O(\seps)$ can be treated as a
perturbation when $\varepsilon$ (suppressed in the latter notations)
is small enough. Further in the paper, the bold typeface marks the
$n+1+m$ dimensional quantities.

If one truncates the normal form Hamiltonian $H_{\rm nf}$ by
dropping the terms $f_{\rm nf}$ and $\g_{\rm nf}$ in the formula
(\ref{nfh}), the action $\iota$ is flow-invariant. For $\iota=0$,
one can separate a natural system of $1+m$ degrees of freedom, whose
Hamiltonian can be written as
\begin{equation}
K(y_0,\ldots y_m)+U(x_0,\ldots x_m),\label{slow}
\end{equation}
where $K(y)$ is a symmetric positive definite quadratic form in
$y\in\R^{1+m}$ and $U(x)$ -- a scalar function on $\T^{1+m}$. In the
sequel, saying that some function is a ``function on a torus''
implies $2\pi$-periodicity of this function in corresponding
variables.

In the simple resonance case $m=0$, one can show that inherent in
the normal form dynamics is the exponentially small separatrix
splitting phenomenon, see \cite{LMS}, \cite{RT}.

In order to show how the exponentially small splitting theory can be
built in the multiple resonance case $m\geq1$, let us consider a
simple model, which generalizes the so-called Thirring model for a
simple resonance, see \cite{Ga}.

Namely, we study the following model Hamiltonian:
\begin{equation}\label{model}H_\mu(\p,\q) = \spr{\omega}{\iota}
+{1\over2}\sum_{j=1}^n\iota_j^2+
H_{1+m}(y_0,\ldots,y_m,x_0,\ldots,x_m) + \mu
V(\varphi,x_0,\ldots,x_m),\end{equation} where $\omega$ is
Diophantine, $V$ is a real-analytic function on $\T^{n+1+m}$ and
$\mu$ is a small parameter. Specifically, for some strictly
increasing sequence of positive reals $l_0,l_1,\ldots l_m,$ let
$H_{1+m}(y,x)$ have the following form:
\begin{equation}\label{um}\begin{array}{lll}
H_{1+m}(y,x)&= &K_{1+m}(y) + U_{1+m}(x), \\ \hfill \\
K_{1+m} (y_0,y_1,\ldots,y_m)&=&{1\over 2}\sum_{i=0}^m {y_i^2\over l_i^2}, \\
\hfill \\
U_{1+m}(x_0,x_1,\ldots,x_m)  &= & \sum_{i=0}^m l_i
\left(\prod_{j=i}^m\cos{x_j}-1\right).
\end{array}
\label{ns}
\end{equation}
Geometrically, the natural system (\ref{ns}) can be visualized as a
``toroidal" pendulum, i.e. a particle of unit mass, confined to move
on the surface of a "vertically standing" in $\R^{2+m}$ torus of
dimension $1+m$, with principal radii $l_m,\ldots,l_0,$ under the
influence of gravity with the free fall acceleration equal to $1$.
Mechanically, the case $m=1$ can be realized as a double pendulum,
whose shorter arm of length $l_0$ is attached to the terminal point
of the longer arm of length $l_1$ and moves in a circle, which rests
upon the longer arm.

On the energy level $H_{1+m}^{-1}(0)$, the origin $O=(0,0)$ is a
single hyperbolic fixed point, with the characteristic exponents
\begin{equation}\lambda_i={1\over l_i}\sqrt{\sum_{j=0}^i l_j},\, i=0,\ldots,m.\label{lms}\end{equation}
Suppose the sequence $l_j$ grows rapidly enough to ensure
\begin{equation}\lambda_0>\max(\lambda_1,\ldots,\lambda_m)\label{lmsin}.\end{equation}
In addition, we have to assume that
\begin{equation}
\inf_{k\in\Z^m_+}\left|\lambda_0-\sum_{j=1}^m k_j
\lambda_j\right|>0, \label{snr}
\end{equation}
where $\Z_+$ denotes non-negative integers.

The origin is connected to itself by a family of homoclinic orbits,
in fact there exist homoclines representing each homotopy class on
$\T^{1+m}$ for the geodesic flow generated by the corresponding
Jacobi metric, degenerate at $x=0$, see \cite{BR}. Some of these
homoclinic orbits, or separatrices, are patently obvious: let
$x_j=y_j=0, \,\forall
j\in\{0,\ldots,m\}\setminus\{i\},\,x_i(t)=4\arctan e^{\pm\lambda_i
t}.$ These orbits correspond to homoclinic geodesics forming the
basis of the fundamental group of the torus $\T^{1+m}$ (modulo the
sign in the exponential which bears witness to the reversibility of
$H_{1+m}$, identifying the upper or lower separatrix branch, where
$y_i$ retains its sign). Consider the orbit with $i=0$, call it
$\gamma$. This orbit leaves and arrives back to the fixed point in
the maximum expansion/contraction direction, corresponding to the
Lyapunov exponent $\lambda_0$. In order to take both of the orbit's
branches into account, let us represent $\gamma$ as follows:
\begin{equation}\gamma =\{x_{1}=\ldots=x_{m}=y_1=\ldots=y_{m}=0,\,y_0=2\sin(x_0/2)
\equiv \psi(x_0), \,x_0\in(0,2\pi)\cup(2\pi,4\pi)\}.
\label{hor}
\end{equation}
Observe that the existence of the two branches of $\gamma$, on each
of which $y_0$ retains its sign, is reflected by
$2\pi$-antiperiodicity of the ``separatrix function'' $\psi$:
$\psi(x_0)=-\psi(x_0+2\pi)$. To reflect this fact, it will be
further convenient to deal with $x_0\in\T_2\equiv\R/4\pi\Z$ rather
than $\T=\R/2\pi\Z$. In particular, addition of values $x$ is
further meant to be $ mod \,(4\pi)$.

Clearly, the orbit $\gamma$ belongs to both the unstable and the
stable $1+m$ dimensional invariant Lagrangian manifolds $W^{u,s}_O$
of the fixed point at the origin. If $m\geq 1,$ the flow of the
Hamiltonian $H_{1+m}$ is non-integrable\footnote{The flow of
$H_{1+m}$ should possess no global analytic first integral other
than the energy, unless $K_{1+m}$ is diagonal and $U_{1+m}$
separated, see \cite{DK}. Transversality of the intersection of the
manifolds $W_O^{u,s}$ along $\gamma$ (to be shown) is in turn an
onset for non-integrability, see \cite{Do}. For the general
variational approach to homoclinic trajectories in natural systems
see \cite{BR}.}.

Global geometry of the manifolds $W^{u,s}_O$ is complicated. Locally
near the origin however, the germs $W^{u,s}_{O, {\rm loc}}$ of the
manifolds $W^{u,s}_{O}$ are diffeomorphic to $m+1$ disks, tangent at
$O$ to the unstable and stable manifolds of the flow, linearized
near the origin.

We shall further show that $\gamma$ arises as a transverse
intersection of the manifolds  $W^{u,s}_O$. Let us call
$W^{u,s}_\gamma$ the localizations of these manifolds in the
neighborhood of $\gamma$. As the orbit $\gamma$ takes off
from/arrives at the fixed point in the maximum expansion/contraction
direction, it itself turns out to be hyperbolic within the manifolds
$W^{u,s}_\gamma$. Indeed, on the ``vertical torus'' in $\R^{2+m}$,
the coordinate directions $x_1,\ldots,x_m$ are the main curvature
directions away from $\gamma$.

Let us further change the notations $(x_0,y_0)$ to $(x,y)$,
$(x_1,\ldots,x_m)$ to $z$ and $(y_1,\ldots,y_m)$ to $\bar{z}$ and
restrict $|z|\leq r_0$ for some $0<r_0<1$. Then
\begin{equation}
H_{1+m} (y,\bar z,x, z) = {y^2\over2l_0^2} + l_0(\cos x-1)
+\sum_{i=1}^m\left[{\bar z_i^2\over 2l_i^2} - \left(l_0\cos{x}
+\sum_{j=1}^i l_j\right) {z_i^2\over 2}  \right]+ O_4(z;x),
\label{rstr}
\end{equation}
The semicolon in the symbol $O_4(z;x_0)$ means that the term in
question is $O(\|z\|^4),$ uniformly in $x_0$, $\|\cdot\|$ standing
further for the Euclidean norm, to be used intermittently with the
sup-norm $|\cdot|$. This notational convention will be used further
on, the parameters following the semicolon often being omitted.

Before formulating the main result, let us give some geometric
description of what we are going to claim. Lifted into the phase
space of the truncated Hamiltonian $H_{\mu}$, where $\mu=0$, the
orbit $\gamma$ gives rise to an isotropic $n+1$ dimensional
invariant manifold, which is topologically a cylinder over the
$n$-torus. Let us denote this cylinder as ${\cal C}_{O}$. Along
${\cal C}_{O}$, there intersects -- degenerately in $n$ directions
corresponding to the rotators' variable $\varphi$ -- a pair of
invariant Lagrangian manifolds ${\cal W}^{u,s}_{O}$, both containing
an invariant whiskered $n$-torus ${\cal T}_{O},$ located at
$(\p,x,z)=(0,0,0).$ On the torus itself, the truncated flow is
quasiperiodic, with the Diophantine frequency $\omega.$ Owing to the
fact that ${\cal C}_{O}$ has two branches, plus the fact that the
trajectories on ${\cal C}_{O}$ are bi-asymptotic to the invariant
torus ${\cal T}_O$, we shall technically refer to ${\cal C}_{O}$ as
a ``bi-infinite bi-cylinder'', yet tending to avoid this rhetoric,
as much as possible.

We study how the presence of the coupling term $V$ in (\ref{model})
is to affect the above described geometric structure and obtain
qualitative estimates for the degeneracy removal effect. As far as
the Hamiltonian $H_{1+m}$ is concerned, the condition (\ref{lmsin})
results in local hyperbolicity of the orbit $\gamma$ within the
manifolds $W^{u,s}_O$ (recall that their localizations near $\gamma$
are denoted as $W^{u,s}_\gamma$). I.e. the germs $W^{u,s}_{O, {\rm
loc}}$ will be contained in the closure of $W^{u,s}_\gamma$ for the
unstable/stable manifolds respectively. Let us denote ${\cal
W}^{u,s}_{\gamma}\cong\T^n\times W^{u,s}_{\gamma},$ the lifting of
the manifolds $W^{u,s}_\gamma$ into the phase space of the truncated
Hamiltonian $H_{\mu},$ when $\mu=0.$ The manifolds ${\cal
W}^{u,s}_{\gamma}$ can be represented by their generating functions
${\cal S}^{u,s}_\gamma(x,z)$ as graphs over the configurations space
variables $(\varphi,x,z)$, where $\varphi\in\T^n$, $|z|<r$ (for some
small enough $r$ to be determined) and
$x\in\T_2\setminus(2\pi-\delta,2\pi+\delta)=[-2\pi+\delta,2\pi-\delta],$
for some positive $\delta<1$.

Then in the perturbed problem we are going to prove the existence of
Lagrangian manifolds ${\cal W}^{u,s}$ representing the analogs of
the manifolds ${\cal W}^{u,s}_{\gamma}$, as far as the Hamiltonian
$H_{\mu}$ is concerned. Moreover, the homoclinic, or bi-infinite
cylinder ${\cal C}_O$ in the truncated system will give rise to a
pair of semi-infinite cylinders ${\cal C}^{u,s}$ in the perturbed
system, each perturbed cylinder containing the invariant whiskered
torus ${\cal T}$, the cylinders ${\cal C}^{u,s}$ themselves being
contained in the Lagrangian manifolds ${\cal W}^{u,s}$ respectively.
The phase trajectories on ${\cal W}^{u}$ will approach ${\cal
C}^{u}$ in negative time; in turn the trajectories
 on ${\cal C}^{u},$ in negative time will approach ${\cal T}$ at a faster rate.
Similar orbit behavior will occur on ${\cal W}^{s}$ and ${\cal
C}^{s}$ in positive time. This is the content of the structural
stability theorem, Theorem \ref{sst} further in the paper.

Moreover, the perturbed Lagrangian manifolds ${\cal W}^{u,s}$ can be
represented as graphs over the configuration space variables, by
adding to the unperturbed generating functions ${\cal
S}^{u,s}_\gamma(x,z)$ respectively, some quantities ${\cal
S}^{u,s}_\mu(\varphi,x,z)$, which are both $O(\mu)$. Then let
 \begin{equation}
 {\cal S}^{u,s}(\varphi,x,z) =  {\cal S}^{u,s}_\gamma(x,z) + {\cal
S}^{u,s}_\mu(\varphi,x,z)\label{ptm}\end{equation}  denote the
generating functions of the perturbed manifolds ${\cal W}^{u,s}$,
respectively.

In order to find these generating functions, we shall describe a
series of canonical transformations, each of which explicitly takes
advantage of the fact that the phase space is a cotangent bundle. In
the sequel, any canonical transformation $\Psi$ will be determined
by some automorphism $\ar$ and closed one-form $dS$ on the base
space. I.e. all the canonical transformations dealt with herein have
the following structure:
\begin{equation}
{\Psi}\,=\,{\Psi}(\ar,S):\ \left\{
\begin{array}{llllllll}
\q&=&\ar(\q'),\\ \p&=&\ti{(d\ar)}\p'+dS(\q),
\end{array}
\right. \label{zed}
\end{equation}
Observe that there is a natural semidirect product structure that on
the pairs $(\ar,S)$, induced by composition.

\medskip\noindent
Hyperbolicity of the orbit $\gamma$ does not suffice to prove
Theorem \ref{sst} however: we also need a special non-resonance (yet
not very restrictive) assumption (\ref{snr}) on the stability
exponents of $H_{1+m}$ at the origin, built into the choice of the
arm lengths $\{l_j\}_{j=0,\ldots,m}$. The latter assumption appears
to be a very special case of the problem of analytic conjugacy
between linearized and non-linear dynamics near a hyperbolic fixed
point, see e.g. \cite{Po}, although for our purposes it suffices
separating the dynamics in a single chosen direction only.

The whiskered torus ${\cal T}_{O}$ and its local unstable and stable
manifold ${\cal W}_{O,{\rm loc}}^{u,s}$ are known to survive small
perturbations without the assumption (\ref{snr}), by the theorem of
Graff, see \cite{Gr}, \cite{Tr}. Namely, in the normal form
(\ref{nfh}) as long as $\omega$ is Diophantine, $U$ possesses a
single non-degenerate absolute maximum, plus the upper left $n\times
n$ minor of the matrix $Q_{\rm nf}$ is nonzero, there exists a
perturbed torus ${\cal T}$ where the flow is conjugate to the linear
flow on its prototype in the truncated system.

We further study the splitting of the unperturbed cylinder ${\cal
C}_{O}$. In order to do so, we introduce the ``splitting function''
\begin{equation}
{\cal D}(\varphi,x,z)= {\cal S}^{u}(\varphi,x,z)- {\cal
S}^{s}(\varphi,x-2\pi,z),\label{spf}\end{equation} which will be
well defined for $\varphi\in\T^n,$ $x\in[-2\pi+\delta,-\delta]\cup
[\delta,2\pi-\delta]$ (recall that addition of $x$ is
$mod\,(4\pi)$), and $|z|<r$. A critical point of ${\cal D}$ would
yield a homoclinic connection to the torus ${\cal T}$, the gradient
$d{\cal D}$ being the ``splitting distance''.

As the manifolds $W^{u,s}_\gamma$ for the Hamiltonian $H_{1+m}$
intersect transversely at $z=0$, the critical points of ${\cal D}$
will lie close to $z=0$, and therefore, the magnitude of the
splitting of the cylinder ${\cal C}_{O}$ can be evaluated in terms
of the derivatives ${D}_{x,\varphi} {\cal D}(\varphi,x,z)$ at $z=0$,
in the properly adjusted coordinate chart\footnote{This is the only
instant in the argument of this paper, where the built into the
model transversality of the intersection of the manifolds
$W^{u,s}_\gamma$ comes into play.}.

\medskip
\noindent The forthcoming Theorem \ref{mth} makes these claims
precise. In order to formulate the theorem, let us introduce some
notation and summarize what analyticity properties are required of
the perturbation $V(\varphi,x,z)$.

For real ${r},\sigma>0$ and $j=1,2,\ldots$ ($j=1$ usually being
omitted) let
\[
\begin{array}{lll}
\B^j_{r} &\=&\{\zeta\in\C^j:\,\|\zeta\|\leq{r}\},\\ \T^j_\sigma & \=
&\{\zeta\in\C^j:\,\Re{\zeta}\in\T^j,\,|\Im{\zeta}|\leq\sigma\}.
\end{array}
\]
For $x\in\T_2$, define a conformal map $s$ and some associated
quantities as follows:
\begin{equation}
s(x)=\int_{\pi}^x{d\zeta\over \psi(\zeta)},
\;\;\;\chi(s)=\psi[x(s)], \;\;\;e=y\chi(s). \label{stime}
\end{equation}
Recall, in the model studied $\psi(x)=2\sin(x/2)$. The map $s(x)$
takes $(0,4\pi)$ to $\R\cup \R+i \pi,$ and the change
$(x,y)\rightarrow(s,e)$ is canonical. The function $x(s)$ is $2\pi
i$-periodic and has singularities at $s=\pm {\pi\over2}i$.

Fix some  $T_0\gg 1$. By construction of the map $s$, for any and
$\rho\in(0,\pi/2)$ any $T\in [T_0/2,T_0],$ the quantities
$x(s),\chi(s)$ are holomorphic functions in the set
$\check\Pi_{\te,\rho}\subset \C/2\pi i$, obtained by throwing out of
$\C$ horizontal rectangles with half-axes
$(2T_0-T)\times(\pi/2-\rho)$, centered at $\pm {\pi\over2}i$.
Namely, let
\begin{equation}
\begin{array}{lll}
\check\Pi_{\te, \rho} &= & \Pi_{\te,\rho}\cup -
\Pi_{\te,\rho},  \\ \hfill \\
\Pi_{\te,\rho}&= &\{ \Re s \leq T,|\Im s|\leq\rho\} \cup  \{ \Re
s \leq T,|\Im s -\pi|\leq\rho\} \cup\{ \Re s\leq T-2T_0\};\;\mbox{ also let}\\ \hfill \\
\hat\Pi_{\te, \rho}&=& \{s\in\C:\,|\Re s|\leq T,|\Im s|\leq\rho\}.
\end{array}
\label{lbd}
\end{equation}
The domains $\Pi_{\te,\rho}$ are further referred to as
semi-infinite bi-strips, their size increasing with $(T,\rho)$, with
$\rho<{\pi\over2}$. Bi-strips $\check\Pi_{\te,\rho}$ are
bi-infinite, while $\hat\Pi_{\te, \rho}$ is simply an
origin-centered horizontal rectangle in $\C$, with semi-axes
$(T,\rho).$

Let
\begin{equation}
{\cal
C}_{\sigma,\te,\rho}=\T^n_\sigma\times\Pi_{\te,\rho},\;\;{\mathfrak
C}_{\sigma,\te,\rho,{r}}={\cal C}_{\sigma, T,\rho}\times \B^m_{r}
\label{cylinders}
\end{equation}
(and in the same fashion $\check{\cal C}_{\sigma,
T,\rho},\,\check{\mathfrak C}_{\sigma,\te,\rho,{r}}$ or $\hat{\cal
C}_{\sigma, T,\rho}, \,\hat{\mathfrak C}_{\sigma,\te,\rho,{r}}$) be
referred to as complex semi-infinite (bi-infinite or finite)
bi-cylinders for ${\cal C}$ and extended bi-cylinders for the
notations ${\mathfrak C}$. In qualitative argument, the analyticity
indices as well as ``bi-'' rhetoric are avoided.

Let us now quote the main assumption.
\begin{assumption} \label{mass} Assume the non-resonance conditions (\ref{dio}) and (\ref{snr}).
Suppose the real-analytic function $V(\varphi,x,z)$ is such that
$V(\varphi,x(s),z)$ is holomorphic and uniformly bounded by $1$ in
$\check{\mathfrak C}_{\sigma_0,T_0,\rho_0,r_0}$ for some initial set
$(\sigma_0,T_0,\rho_0,r_0)$ of analyticity parameters.
\end{assumption}

The main result of the paper is the following theorem.

\begin{theorem}
Under Assumption \ref{mass}, take $T=T_0-1$ and any positive
$\rho<\rho_0,$ $\sigma<\sigma_0,$ let $\delta\sim\log T$. Suppose
\begin{equation}
r< c_1 \min[(\rho_0-\rho), (\sigma_0-\sigma)],\;\;\;\;\;\mu < c_2
[r\vartheta |\omega|^{-1} (\sigma_0-\sigma)^{\tau_n}]^2,
\label{scond}
\end{equation} for some constants $c_{1,2}>0$,
determined by the separatrix function $\psi$ as well as the
quantities $n,\tau_n,m,\sigma_0,T_0,\rho_0,r_0,l_0,\ldots l_m.$
\begin{enumerate}
\item Some level set of $H_{\mu}$, with energy $O(\mu),$
contains an invariant partially hyperbolic $n$-torus ${\cal T}$,
where the flow is conjugate to linear, with the rotation vector
$\omega$. At the torus ${\cal T}$, there intersects a pair of
isotropic manifolds ${\cal C}^u$ and ${\cal C}^s$, which are
contained respectively in the global unstable and stable manifolds
of ${\cal T}$. The manifolds ${\cal C}^u$ and ${\cal C}^s$ are
contained respectively in a pair of Lagrangian manifolds ${\cal
W}^{u,s}$, which are graphs of closed one-forms, with the generating
functions ${\cal S}^{u,s}(\varphi,x,z)$, as in (\ref{ptm})
respectively, such that the quantities ${\cal S}_\mu
^{u,s}(\varphi,x(s),z)$ are holomorphic and uniformly bounded by
$O(\mu)$ for $(\varphi,x,z)\in {\mathfrak C}_{\sigma,\te,\rho,r}$.
The cohomology classes $\xi^{u,s}\in H^1(\T^n,\R)\cong\R^n$ of the
one-forms $d{\cal S}^{u,s}$ are equal to each other.

\item The distance between the manifolds
${\cal W}^{u,s}$ can be measured by the exact one-form $d{\cal D}$,
defined by (\ref{spf}). There exist a coordinate chart
$(\varphi',x',z')\in \T^n\times [\delta,2\pi-\delta] \times\B^m_r$,
obtained by a near-identity change of variables
$(\varphi',x',z')=\ar(\varphi,x,z)$ from the original coordinates in
(\ref{model}), such that in the chart $(\varphi',x',z')$, the
function ${\cal D}$ satisfies the following PDE:
\begin{equation}\label{lhje}
\psi(x'){\partial {\cal D}\over\partial x'} + \spr{\omega}{{\partial
{\cal D}\over\partial\varphi'}}+ \spr{z'}{L'[{\cal D}]}=0,
\end{equation}
where $L'$ is a linear first order differentiation operator.

\item In the above chart, namely for $(\varphi',x',z')\in \hat{\mathfrak
C}_{\sigma, T,\rho,r}$, the function ${\cal D}(\varphi',x',z')$ is
bounded by $O(\mu)$. Let ${\cal D}(\varphi',x',z')={\cal
D}_0(x',\varphi')+O(z').$ Then the quantity  ${\cal
D}_0(x',\varphi')$ can be written as a $2\pi$-periodic function on
$\T^n$ of \begin{equation}\alpha=\varphi' -\omega s(x'), \mbox{ i.e.
} {\cal D}_0(x',\varphi')={\fS}(\alpha).\label{alph}\end{equation}

\item The manifolds ${\cal W}^u$ and ${\cal W}^s$ intersect at least
$2n+2$ orbits, biasymptotic to ${\cal T}$.
\end{enumerate}
\label{mth}
\end{theorem}
Let us show for the moment that the conclusions (iii) and (iv) of
Theorem \ref{mth} are straightforward consequences of (i) and (ii).
Indeed, (\ref{lhje}) implies that ${\cal D}_0(x',\varphi')$ has to
satisfy the linear PDE
$$
\psi(x'){\partial {\cal D}_0\over\partial x'} +
\spr{\omega}{{\partial {\cal D}_0\over\partial\varphi'}}=0,
$$
and $\psi(x'){\partial \over\partial x'}={\partial\over \partial
s}$, where $s(x')$ comes from (\ref{stime}). Then ${\cal
D}_0(x',\varphi')=\fS(\alpha)$ follows, as the form $d{\cal D}_0$ is
exact, i.e. ${\cal D}_0(x',\varphi')$ is $2\pi$-periodic in
$\varphi'\in\T^n$. It follows that the set of the critical points of
the function $\fS(\alpha)$ in the coordinate plane $z'=0$ determines
the trajectories, biasymptotic to the torus ${\cal T}$. The minimum
number  $n+1$  of critical points of $\fS$ is the
Ljusternik-Schnirelmann characteristic of $\T^n$, which equals
$n+1$. It gets doubled in the statement (iv), because one can
restrict $x\in[-2\pi+\delta,-\delta]$, considering the lower
separatrix branch and replicate the statement (ii).

Theorem \ref{mth} has an immediate corollary, implying the estimate
claimed in the Abstract and exponential smallness of the splitting
distance if $\omega\rightarrow\omeps$ for a small $\varepsilon$. We
do not elaborate on various parameter relations in this paper (they
can be all made to depend on $\varepsilon$ in order to approach the
lower bounds' problem for the exponentially small splitting
distance) as the situation here would be the same as it is in the
simple resonance case, regarding which see e.g. \cite{DG} and the
references therein. The willing reader may synthesize these
relations using the fact that the parameter relations in the
forthcoming technical statement Theorem \ref{sstprime} are
supposedly optimal. Going carefully through the latter theorem, one
can derive what precisely the symbols $O(\mu)$ imply regarding the
other parameters of Theorem \ref{mth}.
\begin{corollary}\label{mcr}
For $k\in\Z^n\setminus\{0\}$, the Fourier coefficients ${\fS}_k$ for
of the function ${\fS}(\alpha)$ satisfy the estimate
\begin{equation}
|{\fS}'_k|\,\leq\,O(\mu)  \cdot  e^{-\rho |k\cdot\omega|-
|k|\sigma}. \label{exx}
\end{equation}
If $(\varphi,x,z)$ are the original coordinates and
$(\varphi',x',z')=\ar'(\varphi,x,z)$ is the change of variables,
described by Theorem \ref{mth} (ii), then for all $(\varphi,x,z)\in
\T^n\times[\delta, 2\pi-\delta]\times [-r,r]^m$, one has a uniform
bound
$$
|{\cal
D}_0\circ\ar'(\varphi,x,z)|\,\leq\,O(\mu)\cdot\sum_{k\in\Z^n\setminus\{0\}}e^{-\rho
|k\cdot\omega|- |k|\sigma}.
$$
\end{corollary}
Observe that Theorem \ref{mth} is also valid in the simple resonance
case $m=0$, when the manifolds ${\cal C}$ and ${\cal C}^{u,s}$ are
Lagrangian, rather than isotropic. The simple resonance case has
been exposed in detail in \cite{RT}, see also the references
contained therein. In order to keep the ideas clear and not to rival
the latter reference length-wise, intermediate steps in deriving
many estimates in this paper are omitted, no intermediate
analyticity parameters are explicitly introduced, and the constants
$c_{1,2}$ can get smaller from one statement to another. The reader
is often referred to \cite{RT}, as many aspects of the routine here
mimic the $m=0$ case. The proof of how Corollary \ref{mcr} follows
form Theorem \ref{mth} can also be found in the latter reference, as
well as \cite{Sa}, \cite{LMS}.

\medskip
\noindent {\em Remark.} In the perturbation
$V(\varphi,x_0,\ldots,x_m)$ in (\ref{model}), let
$V_0(\varphi,x_0)=V(\varphi,x_0,0\ldots,0)$, and suppose it vanishes
if $x_0=0$. Then one can formally set up the Melnikov integral
\begin{equation} \mathfrak M(\alpha)=\int_{-\infty}^{+\infty} V_0(\alpha+\omega
t,x_0(t))dt,\;\;\alpha\in\T^n, \label{meln}\end{equation} and ask
whether this quantity is a bona fide $C^2$ approximation for
$\fS(\alpha)$. (A routine calculation shows that setting
$\omega\rightarrow\omeps$ results in the Fourier coefficients
$\mathfrak M_k,$ equal to the right-hand side of the bound
(\ref{exx}), with $\rho={\pi\over2}$.) We do not study this
undoubtedly important issue here (see. e.g. \cite{DG}, \cite{LMS}
for thorough discussion). However, the theory developed further
suggests that there are no extra difficulties arising in this
respect in the multiple versus the simple resonance case. In
particular, the ``easy'' case $n=1$ involving no small divisors
should be similar in this respect to the Thirring model, cf.
\cite{Ga}.

\medskip
\noindent {\em Remark.} Note that contrary to the simple resonance
case, where there exists a large body of literature on exponentially
small bounds for the splitting of separatrices, see e.g. \cite{LMS}
and the references therein, multiple resonances have been usually
approached via the normal forms, alias averaging method. The latter
technique (see \cite{LMS}, \cite{PT}) is not very explicit
geometrically, however as \cite{PT} points out, it does enable one
to obtain exponentially small upper estimates with sharp constants,
which come from dynamical considerations regarding the analyticity
domains, if not to relate these estimates directly to the splitting
of separatrices.

\subsection*{2. Unperturbed system analysis}
\addtocounter{section}{1} \setcounter{equation}{0}

In this section, for the sake of clarity, we confine ourselves to
the case $m=1$ only; the extension to $m>1$ is transparent. Thus, in
this section, let $l_0=1$, $l_1=l>1$. Let us further assume that $l$
is such that $\lambda={\sqrt{1+l}\over l}<1$ and for $k\in\Z_+$
\begin{equation} |k\lambda-1|\geq {\lambda\over 10}, \;\;\forall\,k\in\Z_+,\label{snr1} \end{equation}
a particular case of the condition (\ref{snr}).

The Hamiltonian $H_{1+m}$ given by (\ref{rstr}), for $m=1$ turns
into
\begin{equation}
H_{2} (y ,\bar{z},x,z) = {y^2\over2}+(\cos{x}-1) + {\bar{z}^2\over
2l^2}-(l+\cos{x}) {z^2\over 2} + O_4(z;x). \label{chgonce}
\end{equation}

Let us further change $y \rightarrow \pm\psi(x)+y$, recall that
$\psi(x)=2\sin(x/2)$. Clearly, the choice of the sign as $+$
corresponds to localization near the orbit $\gamma$ as part of the
unstable manifold of the fixed point $O$, while the $-$ sign would
imply doing it near $\gamma$ as part of the stable manifold. Let us
call the resulting Hamiltonians $H_{2,\pm\psi}$ as follows:
\begin{equation}
H_{2,\pm\psi} (y ,\bar{z},x,z) = \pm y
\psi(x)+{y^2\over2}+{\bar{z}^2\over 2l^2}-(l+\cos{x}) {z^2\over 2} +
O_4(z;x). \label{chg}
\end{equation}
Observe that the Hamiltonians $H_{2,\pm\psi}$ have resulted from
$H_2$ after the canonical changes with the generating functions
\begin{equation}
{\cal S}^{\pm}_{\psi}=\pm\int_0^x \psi(\zeta)d\zeta, \label{psis}
\end{equation}
with the $+$ sign for the stable and the $-$ sign for the unstable
manifolds, respectively. Further calculations will be quoted for
mostly $H_{2,+\psi}\equiv H_{2,\psi}$ only.

The Hamiltonian $H_{2,\psi}$ is now a function of
$x\in\T_2=\R/4\pi\Z$, rather than $\T$, with two singular points,
where $x=0,2\pi$. To identify them $H_{2,\psi},$ retains a symmetry:
\begin{equation}\label{spt1}
H_{2,\psi}(y ,\bar z,x,z) = H_{2,\psi}(y
+2\psi(x),\bar{z},x+2\pi,z).
\end{equation}
This symmetry was called  the ``sputnik property'' in \cite{RT} and
was used to validate the analogue of the claim $\xi^u=\xi^s$ of
Theorem \ref{mth} in the case of a simple resonance. Here, as we are
dealing with the specific model, we will not use this property
explicitly to prove this claim, but rather will later compare the
pair $H_{2,\pm\psi}$ (which is anyway tantamount to the same trick,
used in the beginning of section 4).

To identify the unstable/stable manifolds $W^{u,s}_\gamma$ of the
orbit $\gamma$, we will be looking along the $x$-axis at the
unstable/stable manifold of the singular point at $x=0$ for the
Hamiltonians $H_{2,\pm\psi}$, respectively. This proves convenient
and suggests that in general the regular description of the
manifolds $W^{u,s}_\gamma$ is likely to fail in the neighborhood of
$x=2\pi$.

Let us linearize the flow of the Hamiltonian $H_{2,\psi}$ near the
orbit $\gamma$, whereupon $x(t)\equiv x_0(t)=\pm 4\arctan e^{ t}$.
For the infinitesimal increments $(\hat x, \hat y , \hat z,
\hat{\bar{z}}),$ one gets the system of equations
\begin{eqnarray} \label{lin1} \dot{\hat x}\;=\; D\psi[x_0(t)]\hat x+\hat
y , & &\dot{\hat y}\;=\;-D\psi[x_0 (t)]\hat y , \\ \nonumber \\
\dot{\hat z} \;=\;\hat{\bar{z}}/l^2, & & \dot{{\hat{\bar z}}} \;=
\;(l-1+\tanh^2 t)\hat z.\label{lin2}
\end{eqnarray}
The tangent space to  $W^{u}_\gamma$  at the points on $\gamma$ will
be spanned by the vectors -- solutions of the latter system of
linear ODEs, which vanish as $t\rightarrow -\infty$.

The two pairs of equations (\ref{lin1}), (\ref{lin2}) are uncoupled.
As far as (\ref{lin1}) is concerned, there is an obvious solution
$\hat{x}(t)=\dot{x}_0(t)\sim 1/\cosh t,\,\hat{y}(t)=0$, which
vanishes at both $t\rightarrow\mp\infty$. I.e. one tangent direction
to $W^{u}_\gamma$ at a point $(x,z,y ,\bar{z})=(x,0,0,0)$ is always
$(1,0,0,0)$, in the direction collinear with $\gamma$ itself.

Equations (\ref{lin2}) will clearly have no solutions vanishing at
both $t\pm\infty$, as the coefficients therein retain their sign for
all $t$. However the system certainly does have a solution
$(\hat{z}^u(t),{\hat{\bar z}}{}^u(t))$, defined for $t\leq T_0$ for
some $T_0\gg 1$, which as $t\rightarrow-\infty$ approaches the
trivial $(\hat{z},\hat{\bar{z}})\equiv(0,0)$ (as well as another
solution, which is defined for $t\geq-T_0$ and vanishes at
$t\rightarrow+\infty$).

To construct the unstable solution $(\hat{z}^u(t),{\hat{\bar
z}}{}^u(t))$, one may set ${d\over dt}=\psi(x){d\over dx},$
$1-\tanh^2 {t}=\cos{x}$ and construct the germ of the solution in
question locally as a Taylor series in $x$ near $x=0$;  by linearity
of the equations (\ref{lin2}) and boundedness of their coefficients,
the continuation of these germs over a finite time interval does not
itself pose any problem\footnote{As a matter of fact, equations
(\ref{lin2}) represent a second order linear ODE of generalized
Legendre type, whose general solution can be found explicitly in
terms of the associated Legendre functions.}.

Observe that given $(\hat{z}^u(t),{\hat{\bar z}}{}^u(t)),$ one can
let $(\hat{z}^s(t),{\hat{\bar z}}{}^s(t))=
(\hat{z}^u(-t)-2\pi,-{\hat{\bar z}}{}^u(-t))$ for the result of the
similar procedure with respect to the Hamiltonian $H_{2,-\psi}$. For
no $t\in[-T_0,T_0]$ can the  vectors $(\hat{z}^u(t),{\hat{\bar
z}}{}^u(t))$ and $(\hat{z}^s(t)+2\pi,{\hat{\bar z}}{}^s(t))$ be
parallel, or there would exist a solution of (\ref{lin2}),
biasymptotic to zero. For the Hamiltonian $H_2$ from
(\ref{chgonce}), the existence and transversality of the
intersection along the orbit $\gamma$ of a pair of manifolds
$W^{u,s}_\gamma$ (defined in the neighborhood of $\gamma$)
essentially follow. A quantitative statement of this fact is to be
given shortly. So far observe by comparing the coefficients in the
linear equations (\ref{lin2}) retain the sign, both vectors
$(\hat{z}^u(t),{\hat{\bar z}}{}^u(t))$ and $(\hat{z}^s(t),{\hat{\bar
z}}{}^s(t))$ in the $(z,\bar{z})$ plane never have a slope too close
to horizontal. More precisely, their slope in absolute value will be
contained in the interval $[l\sqrt{l-1},l\sqrt{l+1}]$, which can be
seen from (\ref{til}) below.

Let us show how the manifold $W^u_\gamma$ can be constructed, the
analysis for $W^s_\gamma$ gets modified in the obvious way. Make a
change $\bar{z}\rightarrow \bar{z} + \lambda_u(x)z$, where
$\lambda_u(x)$ determines the direction of the solution vector,
vanishing at $x=0$ (i.e. $t\rightarrow-\infty$). The quantity
$\lambda_u(x)\in [l\sqrt{l-1},l\sqrt{l+1}]$ is well defined for
$x\in\T_2\setminus (2\pi-\delta,2\pi+\delta)$ for some $0<\delta<1$,
where $\delta \approx \ln T_0$.

The change $\bar{z}\rightarrow \bar{z} + \lambda_u(x)z$ is not
canonical, to make up for it one also has to change $y \rightarrow y
+{1\over2}{d\lambda_u(x)\over d x} z^2$. In other words, this is a
canonical change with the generating function
\begin{equation}
{\cal S}_{\gamma,0}^u (x,z)= {1\over 2}\lambda_u(x) z^2.\label{lins}
\end{equation}
Then the Hamiltonian $H_{2,\psi}$ in (\ref{chg}) transforms to
\begin{equation}H_{2,u}(x,y)=y \psi(x)+{y^2+l^{-2}\bar{z}^2\over 2}
+{1\over 2} \tilde{\lambda}_u(x) z^2
 + l^{-2}\lambda_u(x)z \bar{z} + yO_{2}(z;x) +
 O_4(z;x),\label{linham}\end{equation}
where
\begin{equation}\label{til}\tilde{\lambda}(x)\;=\;\psi(x){d\lambda_u(x)\over
dx}-(l+\cos x)+l^{-2}\lambda_u^2(x). \end{equation} It follows by
construction -- or directly from (\ref{lin2}) -- that
$\tilde{\lambda}_u(x)\equiv0$. Thus the quantity $\Lambda_u(x)\equiv
l^{-2}\lambda_u(x)$ multiplying $z \bar{z}$ in (\ref{linham}) is
always positive, never exceeding $\lambda=\Lambda_u(0)$; recall that
$\lambda<1$, also cf. (\ref{snr1}).

Finally, the last two terms in (\ref{linham}) can be regarded as a
perturbation, provided that ${r}$ is small enough.

The phase space of the Hamiltonian $H_{2,u}$ is $T^*((\T_2\setminus
(2\pi-\delta,2\pi+\delta))\times[-r,r])$. If there were no two last
terms in (\ref{linham}), the  manifold $W^u_\gamma$ would be given
by the zero section $(y,z)=(0,0)$ of the bundle. However, for small
$r$, the last two terms in (\ref{linham}) can be regarded as a
perturbation and dispensed with, owing to the following lemma.

\begin{lemma} \label{straight}
Given $\rho<\rho_0$, $T\leq T_0-1$, there exists a constant $c_1>0$,
depending only on the parameter set $(\rho_0,T_0,r_0)$ and
$\lambda$, such that for $r< c_1(\rho_0-\rho)$, there exists some
reals $\delta,\kappa=O(1)$ in $(0,1)$ and a canonical near-identity
transformation $\Psi_r^{u}$, such that the Hamiltonian
(\ref{linham}) can be cast into the following normal form:
\begin{equation}
H_{\gamma,u}(y,\bar{z},x,z)=y \psi(x)+ \Lambda_u(x)z \bar{z} +
 O_2(y,\bar{z}),\label{unpham}\end{equation}
 valid for $|y|,|\bar{z}|\leq\kappa$, $|z|\leq r$ and $x$ such that
 $\Re{x}\in [-2\pi+\delta,2\pi-\delta]$ and $s(x)\in \Pi_{T,\rho}$.

The transformation $\Psi_r^{u}$, for $p=(y,\bar z)$ and
 $q=(x,z)$ can be written in the following form:
 \begin{equation}
\Psi_r^{u}= \Psi_r^{u}(b_r^{u},{\cal
S}_{r}^{u}):\;\left\{\begin{array}{llccc}
q&=&q'&+&b_r^{u}(q'),\\
p&=&\ti{[{\rm id}+db_r^{u}(q')]}p' &+& d{\cal
S}_r^u(q),\end{array}\right.\label{strtr}\end{equation} where the
quantities $b_r^{u}(x,z)$ and ${\cal S}_{r}^{u}(x,z)$ are both
$O_2(|x|+|z|)$.
\end{lemma}
We omit the proof of the Lemma, as it follows as a particular case
of the forthcoming Theorem \ref{sst}, in the
 case when there are no $\varphi$-dependencies. The smallness
 condition $r< c_1(\rho_0-\rho)$ follows after routine, but careful examination
 of the proof of Theorem \ref{sst}, see also the quantitative
 estimates in Theorem \ref{sstprime}. As a matter of fact, Lemma \ref{straight} still holds if the term
 $O_4(z;x)$ in (\ref{linham}) gets replaced by $O_3(z;x)$.

Observe that repeating the argument for the Hamiltonian
$H_{2,-\psi}$, the latter would be cast in the following form:

\begin{equation}
H_{\gamma,s}(y,\bar{z},x,z)= - y \psi(x) - \Lambda_s(x)z \bar{z} +
 O_2(y,\bar{z}),\label{unphams}\end{equation}
 where $\Lambda_s(x)>0$ and equals $\lambda$ at $x=0$.
In order to get (\ref{unphams}), the analog of Lemma \ref{straight}
would be preceded by a canonical transformation with the generating
function ${\cal S}^s_{\gamma,0}(x,z)={1\over2}\lambda_s(x)z^2$, cf.
(\ref{lins}).

\renewcommand{\theproposition}{2.\arabic{proposition}}
\addtocounter{proposition}{1} Let us summarize the results of the
analysis in this section by the following proposition
\begin{proposition} For $r$ small enough, the unstable/stable manifolds
$W_\gamma^{u,s}$ for the Hamiltonian (\ref{chgonce}) can be
represented as graphs over the variables $(x,z)$, via the generating
functions
\begin{equation}
{\cal S}_\gamma^{u,s} = {\cal S}_\psi^{+,-} + {\cal
S}_{\gamma,0}^{u,s} +  {\cal S}_{r}^{u,s},
\label{sgam}\end{equation} respectively, the representation being
valid for $\Re{x}\in\T_2\setminus(2\pi-\delta,2\pi+\delta)$,
$s(x)\in \Pi_{\rho,T}$, $\delta\sim \log T$  and $z\in\B^m_r$. Both
${\cal S}_\gamma^{u,s}$ vanish to the second order at $(x,z)=(0,0)$.
The intersection of the manifolds $W_\gamma^{u,s}$ along the orbit
$\gamma$ is transverse, for
$x\in[\delta,2\pi-\delta]\cup[2\pi+\delta, 2\pi-\delta]$.

The Hamiltonian (\ref{chgonce}) can be cast into the forms
(\ref{unpham}), (\ref{unphams}) via canonical changes
$\Psi_\gamma^{u,s}$ respectively, where $\Psi_\gamma^{u,s}=
\Psi_\gamma^{u,s}(a^{u,s}_r, {\cal S}_\gamma^{u,s}),$ with the
near-identity diffeomorphisms $a^{u,s}_r={\rm id}+b^{u,s}_r$ such
that $b^{u,s}_r$ both vanish to the second order at $(x,z)=(0,0)$.
\label{unp}\end{proposition}

The next and most important step towards proving Theorem \ref{mth}
is to get the generating functions ${\cal S}_\mu^{u,s}$ in
(\ref{ptm}) by developing the structural stability theory for a
class of Hamiltonians, which would include
\begin{equation}
H_{u,s}=\spr{\omega}{\iota}+{1\over2}\sum_{j=1}^n \iota_j^2 +
H_{\gamma,\cdot} + V_{u,s}(\varphi,x,z), \label{input}\end{equation}
with $\cdot=u,s$, respectively, and the perturbations $V_{u,s}$
satisfy Assumption \ref{mass}. We can also assume that we can
estimate the partial derivatives of ${D}_{x,z}V_{u,s}$ at
$(x,z)=(0,0)$ in terms of the norm of $V_{u,s}$, as $V_{u,s}$ itself
should be analytic for $|x|\leq\delta$ and $|z|\leq r_0$.

\subsection*{3. Structural stability theory}
\addtocounter{section}{1} \setcounter{equation}{0}

The classical structural stability theorem in problems with small
divisors is due to Kolmogorov \cite{Ko}, alias KAM. It establishes
stability of geometric objects -- invariant Lagrangian tori with
quasi-periodic flow thereupon -- with respect to perturbations of
non-degenerate Hamiltonians localized in the neighborhoods of these
geometric objects. If one is after other geometric objects, say
whiskered tori of lower dimension, with quasi-periodic flow, there
is a theorem of Graff, \cite{Gr}. It was proposed in \cite{RT} that
a proper geometric object to look at in order to set up the
splitting problem near a simple resonance is a semi-infinite
cylinder over a torus. In this section the structural stability
theorem from \cite{RT} is given extra development. Namely, we deal
with "extended cylinders", see (\ref{cylinders}), which appear to be
the proper geometric objects to study in order to describe the
separatrix splitting at multiple resonances. Naturally, all the
above mentioned theorems, as well as the forthcoming Theorem
\ref{sst}, which is the main result in this section, are in the
realm of the abstract implicit function theorem framework of Zehnder
\cite{Z1}, \cite{Z2}.

As this section is the most technical one, the notations in it
(hence also in section 5 and the Appendix) are largely
self-contained. We study the following Hamiltonian, with
$(\p,\q)=(\iota,y ,\bar{z},\varphi,x,z)$ (to justify the earlier
made claim that it suffices to consider $m=1$, let us now take
$z=(z_1,\ldots,z_m),$ for any $m$):
\begin{equation}
H_\omega(\p,\q) = \lambda_0 \psi(x)y + \spr{\omega}{\iota}  +\spr{z}{\Lambda(x) \bar{z}} + O_2(\p;\q), \label{snf}
\end{equation}
under the following basic assumptions:
\begin{enumerate}
\item $\omega$ is Diophantine, for all $\varphi$, the matrix ${D}^2_{\iota\iota}H_\omega(0,0,0,\varphi,0,0)$ is
non-degenerate.
\item For all $x\in \T_2\setminus(2\pi-\delta,2\pi+\delta),$
the real parts of the eigenvalues of the diagonalizable matrix
$\Lambda(x)$ lie in the interval $(0,\lambda_0)$; moreover $\Lambda(0)=
diag(\lambda_1,\ldots,\lambda_m)$ and the condition (\ref{snr}) is satisfied.
\end{enumerate}
Clearly, the simultaneous change $(\lambda_0,\Lambda)$ to
$(-\lambda_0,-\Lambda)$ is not going to violate the principal
conclusions of this section.

\subsubsection*{Notation} Technically it proves convenient to deal
with the  non-compact ``energy-time''
coordinates $(e,s),$ introduced by (\ref{stime}), rather than the
coordinates $(y ,x)$; some notation and formalism are being set up
further.

Let $\fB_\sigma(\T^j)$ be the Banach space of bounded
$2\pi$-periodic scalar functions in each variable, real analytic in
$\T^j_\sigma$, with the sup-norm.

Let $x\in\T_2$. Rather than dealing with a fixed $\psi(x)=2\sin(x/2)$, let us introduce it axiomatically, as
a real-analytic function,
such that $\psi(0)=0,\,D\psi(0)=1.$ Suppose
$\psi(x+2\pi)=-\psi(x)$ and $\psi(x)$ has no other zeroes on the
real line, but integer multiples of $2\pi$.

Define a conformal map $s(x)$ via (\ref{stime}). The map $s(x)$
takes $(0,4\pi)$ to $\R\cup \R+i \pi$ and the change
$(x,y)\rightarrow(s,e)$ is canonical.

By construction of the map $s$, there exists some $T_\psi\gg 1$ and
$\rho\in(0,\pi/2)$ such that for any $T\in [T_\psi/2,T_\psi]$ the
quantities $x(s),\chi(s)$ are holomorphic functions in the set
$\check\Pi_{\te,\rho}\subset \C/2\pi i$, defined by (\ref{lbd}). In
addition, and this is possible by the properties of $\psi(x)$, let
us suppose that $\rho$ is such that for any $s\in \check\Pi_{\te,
\rho}$, there exists a pair of constants $c_\psi$, $C_\psi$ such
that
\begin{equation}
c_\psi e^{-|s|} \leq |\chi(s)|\leq C_\psi
e^{-|s|}.\label{chibounds}\end{equation}

\medskip
\noindent To deal with Hamiltonian functions in $T^*{\mathfrak C}$ let us
introduce some function spaces. For more details, see \cite{RT}.

The function spaces will be characterized in terms of the
analyticity parameters, accommodated into parameter vectors $\fp$ as
follows. Let $\mathfrak{p}=(\sigma,T,\rho)\in\R^3_{++}$. Introduce
partial order $\mathfrak{p}'=(\sigma',T',\rho')\,\leq\,\mathfrak{p}$
if $\sigma'\leq\sigma,T'\leq T,\rho'\leq \rho$. If
$\mathfrak{p}'\leq\mathfrak{p}$ and $|\mathfrak{p}-\mathfrak{p}'|
\equiv\inf(\sigma-\sigma',T-T',\rho-\rho')>0$, write
$\mathfrak{p}'<\mathfrak{p}$. Addition of parameter vectors, as well
as multiplication by positive reals is defined component-wise, as
well as the difference $\mathfrak{p}-\mathfrak{p}'$ for
$\mathfrak{p}'<\fp$. For $\de\in(0,|\fp|),$ the notation
$\fp'=\fp-\de$ means that $\de$ has been subtracted from each
component of $\fp$. In the sequel the components and dimension of
the parameter vectors $\mathfrak{p}$ may vary; more often than ever
we will have $\mathfrak{p}=(\sigma,T,\rho,{r})\in\R^4_{++}$.

Let $\bsj$ be the Banach space -- with the sup-norm -- of bounded
holomorphic functions $u$ on ${\cal C}_{\sigma,T,\rho}$, such that
$u(\varphi,s)=u[\varphi,s(x)]=\tilde u(\varphi,x),$ where $\tilde u$
is bounded and holomorphic in the Cartesian product of $\T_\sigma$
and the pre-image of the strip $\Pi_{T,\rho}$ under the map $s(x)$,
and $\tilde u(\varphi,x)$ vanishes to the $j$th order at $x=0$ (the
index $j=0$ being omitted). Component-wise sup-norm $|\cdot|_{\fp}$
or the equivalent Euclidean norm $\|\cdot\|_{\fp}$ is used for
vector functions.

If $u(\varphi,s)\in \bsj$, a multiplier $\chi^j(s)$ can be factored
out, i.e.
\begin{equation}
u(\varphi,s)=\chi^j(s) v(\varphi,s),\;\;\;v\in \bs,
\;\;\;|v|_{\fp}\,\approx\,|u|_{\fp}, \label{dcp}
\end{equation}
with constants depending on the fixed quantity $\chi$ being
henceforth absorbed into the symbols $\lesssim,\,\approx$, see also
(\ref{chibounds}). For $u\in \bs$, there exists a unique
decomposition
\begin{equation}
u(\varphi,s)\,=\,u_0(\varphi)+u_{1}(\varphi,s),
\hspace{3mm}\mbox{where}\hspace{3mm}u_0\in{\fB}_\sigma(\T^n),\;\;u_{1}\in
{\fB}^1_{\fp}({\cal C}). \label{dcmp}
\end{equation}
Using it, define the average $\langle u \rangle$ ``at infinity'' as
\begin{equation}
\langle u \rangle\,\=\,\int_{\T^n} u_0(\varphi)d\varphi.
\label{aver}
\end{equation}
For $u\in\bsone$, there is an estimate:
\begin{equation}|u(\varphi,s)|\lesssim e^{s}|u|_{\fp}.\label{espdec}
\end{equation}

Let us also introduce a function space $\bsm\cong \bst\times\bs$ of
functions unbounded at infinity as follows: \begin{equation}
u(\varphi,s)\in {\fB}^\wedge_{\fp}({\cal C})\;\mbox{ iff }\;
u(\varphi,s)={v(\varphi,s)\over\chi(s)},\;\; v(\varphi,s)\in \bs.
\label{nbd}\end{equation} The norm on $\bsm$ is defined as
$|v|_{\fp}$. By (\ref{dcmp}) and (\ref{nbd}), for $u\in\bsm$ there
is a decomposition ${\displaystyle
u(\varphi,s)=v_0(\varphi)/\chi(s)+v_1(\varphi,s),}$ for some
$v_1\in\bs$. Also let $\bsg=[\bs]^n\times\bsm$. An element of this
space describes a vector field on $\cyl$ as well as a map
$a(\varphi,s)$ of ${\cal C}_{\fp}$ into ${\cal C}_{\fp'}$, (with
$\fp<\fp'$for the map to be well defined). Namely, if $g\in\bsg$ is
a vector field and $a$ is such a map, then the "new" vector field $d
a^{-1} g\circ a$ is in $\bsgprime$, see \cite{RT}. It is legitimate
to use the Cauchy formula to estimate partial derivatives of $u\in
{\fB}^{\wedge}_{\fp}({\cal C})$, i.e. ${\displaystyle
|du|_{\fp'}\,\lesssim\,\de^{-1}|u|_{\fp},}$ where $\de=\fp-\fp'$.

All the above notations extend in an obvious way to functions on
$\cyle$, by adding a component ${r}$ to the parameter vector $\fp$,
and considering absolutely convergent Taylor series in $z\in
\B^m_{r}$ with the coefficients in the corresponding spaces of
functions on $\cyl$. For instance $\bse$ (with
$\fp=(\sigma,T,\rho,{r})$) becomes an extension of the space $\bs$
(with $\fp=(\sigma,T,\rho)$), and the notation $\bseg$ extends
$\bsg$. The norm in the extended spaces, such as $\bse$ is the sum
of the Taylor series in $z$, where the moduli suprema have been
taken for all the coefficients. The quantity $r$ will not appear
explicitly in the estimates in this section, getting absorbed in the
$\lesssim$ symbols, e.g. for $u\in\bse$, we have $|D_z
u|_{z=0}\lesssim |u|_{\fp}$. To bring it to terms with the fact that
$r$ in Theorem \ref{mth} is actually quite small, cf. (\ref{scond}),
observe that $r$ would further come into play only when the
functions' derivatives are evaluated at $z=0$ via the Cauchy
inequality. But for the functions Theorem \ref{mth} is dealing with,
these derivatives can be estimated in terms of the quantity $r_0$,
which is $O(1)$. The same should be said about the parameters
$(\sigma,T,\rho)$ which are all supposed to be independent of the
parameter characterizing the perturbation size.

The notation $\avr{u}$ for $u\in\bse$ implies that $z$ has been set
to zero, cf. (\ref{aver}). Thus, for $u\in\bseone$, there is a
uniform estimate, cf. (\ref{espdec}):
\begin{equation}|u(\varphi,s,z)|\lesssim e^{s}|u|_{\fp}.\label{espvan}
\end{equation}

Hamiltonian functions on $T^*\cyle$ are given by absolutely
convergent Taylor series with coefficients in $\bse$, in
$\tilde{\p}=(\iota,y,\bar{z})=(\iota,e/\chi(s),\bar{z})$, inside a
complex ball $\B^{n+1+m}_\kappa$. Notation-wise $(e,s)=(0,-\infty)$
corresponds to $y=0$. Let $\bsek$ be the space of such Hamiltonians,
the norm being the sum of the Taylor series in $\tilde{\p}$, where
the norms have been taken for all the coefficients.

\subsubsection*{Structural stability theorem}
What follows is a non-technical formulation of the theorem to keep
its content transparent.
\begin{theorem}\label{sst}
Consider the following Hamiltonian $H_\omega\in\bsek$:
\begin{equation}
H_\omega(\iota,e,\bar{z} ;\varphi,s,z) =  const. + \lambda_0 e +
\spr{z}{\Lambda(s)\bar{z}} + \spr{\omega}{\iota} +
O_2(\tilde{\p};\q), \label{tmp}
\end{equation}
with $\tilde{\p}=(\iota,{y/\chi(s)},\bar{z})$. Assume the following:
\begin{enumerate} \item $\omega\in\R^n$ is Diophantine and the matrix $D^2_{\iota\iota}
H_\omega(0,0,0;\varphi,-\infty,0),\,\forall \varphi$ is
non-degenerate.
\item $\lambda_0>0$, the real parts of all the eigenvalues of
$\Lambda(s)\in\bs,\,\forall s$ lie in the interval $(0,\lambda_0)$;
\item in the decomposition $\Lambda(s)=\Lambda_0+\Lambda_1(s)$, with
$\Lambda_1(s)\in\bsone$, one has
$\Lambda_0=diag(\lambda_1,\ldots,\lambda_m),$ with
$0<\Re\lambda_j<\lambda_0, \,\forall j=1,\ldots,m,$ and the
condition (\ref{snr}) is satisfied by
$\{\lambda_0,\ldots,\lambda_m\}$.
\end{enumerate}
Then $H_\omega$ is structurally stable,  via a canonical
transformation
\begin{equation}
\Psi  \,=\,\Psi   (\ar ,S ): \;\left\{
\begin{array}{llllllll}
\q&=&\ar(\q'),\\ \p&=&\ti{(d\ar )}\p'+dS ,
\end{array}
\right. \label{zed1}
\end{equation}
and for any $\mathfrak{p}'<\mathfrak{p}$, the transformation
$\ar={\rm id}+\br$, with $\br\in\bsegprime$. The one-form $dS$ is
defined by the generating function $S(\q)=\spr{\xi }{\varphi}+{S_0}
(\varphi,s)$, with $\xi \in\R^n,$ ${S_0}\in\bseprime$.
\end{theorem}

Let
\begin{equation}
\begin{array}{ccc}
H&=&H_\omega+V,\\ V(\p,\q) &=& f(\q) + \spr{\g(\q)}{{\p}}
\end{array}\label{pt}
\end{equation}
be a small perturbation of the Hamiltonian (\ref{tmp}). In the
perturbation, suppose $f\in\bse$ and $\g\in\bseg$. How small the
perturbation should be is stated in the forthcoming technical
version of Theorem \ref{sst}, Theorem \ref{sstprime}.

\medskip
\noindent{\em Remark.} An important consequence of the analytic
set-up to be used further is local uniqueness. I.e. given the pair
$(H_\omega,V)$, the pair $(\ar,S)$ defining the conjugacy
transformation $\Psi$ in (\ref{zed1}) is unique.
\medskip\noindent

\medskip
\noindent The proof of Theorem \ref{sst} is given in section 5. Let
us now discuss some implications of the theorem, in the coordinates
$(\varphi,x,z)$, cf. (\ref{stime}), where $\varphi\in\T^n$,
$x\in\T_2\setminus(2\pi-\delta,2\pi+\delta)$, $|z|\leq r$. Let
$\tilde{S},{\tilde{S_0}}$ be the expressions for the generating
functions $\tilde{S},{\tilde{S_0}}$ from Theorem \ref{sst} in these
coordinates.

\begin{corollary}\label{gf}
The Hamiltonian $H$ in (\ref{pt}), as a function of
$(\iota,y,\bar{z},\varphi,x,z)$, on some energy level, possesses an
invariant Lagrangian manifold $\fC$, given by the graph of the
closed one-form $d\tilde S$, where
\begin{equation}\label{genfun}S(\varphi,x,z)=\spr{\xi}{\varphi}+\tilde{S_0}(\varphi,x,z),\end{equation}
and $\tilde{{S_0}}$ is $2\pi$-periodic in the variable $\varphi$.
The manifold $\fC$ contains a partially hyperbolic invariant torus
${\cal T}$, which in turn is contained in an invariant cylinder
${\cal C}\cong \T^n\times [-2\pi+\delta,2\pi-\delta]$.

If the perturbation $(f(\varphi,x(s),z)), {\g}(\varphi,x(s),z))$ in
(\ref{pt}) is such that $f=O_2(|x|+|z|)$ and $\g=O_1(|x|+|z|)$, then
$\xi=0$ and the energy value on the manifold $\fC$ coincides with
the value of the unperturbed Hamiltonian $H_\omega$ thereon.
\end{corollary} Indeed, the first
claim follows from (\ref{tmp}) and (\ref{zed1}) by setting in the
latter formula $\p'=0$. Furthermore, if
$\tilde{S}(\varphi,x(s),z)=S(\varphi,s,z),$ where the latter comes
from Theorem \ref{sst}, then the invariant cylinder ${\cal C}$
arises by letting $z'=0$ in $S'(\varphi',s',z')=S\circ
\ar^{-1}(\varphi',s',z')$, where the transformation $\ar$ also comes
from Theorem \ref{sst}. The torus ${\cal T}$ arises by further
setting $s'=-\infty$.

The second claim follows by observing that a special perturbation,
as described in the Corollary, does not affect the invariant torus
at $(\p,x,z)=0$ and local uniqueness. Alternatively, one can verify
this claim by carefully inspecting the proof of Theorem \ref{sst}.

\subsection*{4. Conclusion of the proof of Theorem \ref{mth}}
\addtocounter{section}{1} \setcounter{equation}{0}

Combining the claims of Proposition \ref{unp}, Theorem \ref{sst} as
well as Corollary \ref{gf} applied to the Hamiltonians $H_{u,s}$ in
(\ref{input}), one immediately establishes the claim (i) of Theorem
\ref{mth}, but for the fact that the invariant manifolds ${\cal
W}^{u,s}$ lie on the same energy level and the fact of equality of
the cohomology class representatives $\xi^{u,s}\in\R^n \cong
H^1(\T^n ,\R)$. Both facts however easily follow by observing that
all the generating functions in (\ref{sgam}) vanish to the second
order at $(x,z)=(0,0)$, where the unperturbed invariant torus is
located, so one can use Corollary \ref{gf}.

Namely, let $\ar^{u,s}_\gamma$ extend the diffeomorphisms
$a^{u,s}_\gamma$ in Proposition \ref{unp}, acting as the identity on
the $\varphi$-variables; let $\Psi^{u,s}_\gamma({\ar}^{u,s}_\gamma,
{\cal S}^{u,s}_\gamma)$, be the corresponding canonical
transformations. Let $\Psi^{u,s}_\mu({\ar}^{u,s}_\mu, {\cal
S}^{u,s}_\mu)$ be supplied by Theorem \ref{sst}, being applied to
the Hamiltonians (\ref{input}), where the quantities
${\ar}^{u,s}_\mu, {\cal S}^{u,s}_\mu$ are viewed as the functions of
$(\varphi,x,z)$ rather than $(\varphi,s,z)$. Let $H_\omega^{u,s}$ be
the results of conjugacy:
$$ H_\omega^{u,s} =
H_\mu\circ\Psi^{u,s}_\gamma\circ\Psi^{u,s}_\mu,$$ respectively for
the unstable and the stable manifolds.

Consider the Hamiltonian $$H'=H_\omega^{u}\circ
(\Psi^u_\gamma)^{-1}\circ \Psi^s_\gamma.$$ By the properties of the
pair $(a_r^{u,s},{\cal S}^{u,s}_\gamma)$ described by Proposition
\ref{unp}, it follows that $H' = H'_\omega+V'$, where $V'=(f',\g')$
is such that $f$ vanishes to the second order and $\g'$ to the first
order at $(x,z)=(0,0)$, while $H'_\omega$ can be regarded as
unperturbed Hamiltonian, in the sense of Theorem \ref{sst}. This
implies that by Theorem \ref{sst} and Corollary \ref{gf}, there
exists a transformation $\Psi'(\ar',S')$, which nullifies the
perturbation $V'$, and the one-form $dS'$ is exact, i.e the
corresponding $\xi'=0$.

Thus
$$
H'\circ\Psi' = (H_\mu \circ \Psi_\gamma^s)\circ
[(\Psi_\gamma^s)^{-1} \circ\Psi^{u}_\gamma\circ\Psi^{u}_\mu
\circ(\Psi^u_\gamma)^{-1} \circ \Psi^s_\gamma\circ \Psi'],
$$
i.e., by uniqueness, the application of Theorem \ref{sst} to the
``stable manifold'' Hamiltonian $H_s=H_\mu \circ \Psi_\gamma^s$ is
effected via the canonical transformation $$(\Psi_\gamma^s)^{-1}
\circ\Psi^{u}_\gamma\circ\Psi^{u}_\mu \circ (\Psi^u_\gamma)^{-1}
\circ \Psi^s_\gamma\circ \Psi'.$$ This transformation is still of
the form (\ref{zed}). Besides, the corresponding generating function
will contain a single ``non-exact'' term $\spr{\xi^u}{\varphi}$,
supplied by $\Psi^{u}_\mu$, as (in the sense of the template
(\ref{zed})) the rest of the transformations in the above chain are
effected by exact one-forms. This proves the claim (i) of Theorem
\ref{mth}.

To prove the claim (ii) of the theorem, substitute $\p=d{\cal
S}^u(\q)$ and  $\p=d{\cal S}^s(\q)$ into the Hamiltonian
(\ref{model}), subtract the result of the latter substitution from
the result of the former one. After substraction has been done, all
the momentum-independent terms are gone, and introducing the
splitting function ${\cal D}$ as in (\ref{spf}), we arrive in the
relation
\begin{equation}
[\psi(x)+O(\mu)]{\partial{\cal D}\over\partial x} +
\spr{\omega+O(\mu)}{{\partial{\cal D}\over\partial \varphi}} +
O(\mu){\partial{\cal D}\over\partial z} + \spr{z}{L[{\cal D}]} = 0,
\label{pthje}\end{equation} where the quantities $O(\mu)$ as well as
the coefficients of the first order linear differential operator $L$
depend on ${\cal S}^{u,s}$. To prove the claim now, it suffices to
solve the vector field conjugacy problem, which ensures the
structural stability of the vector field $\x_0=\left({\partial\over
\partial s}, \spr{\omega}{{\partial\over
\partial \varphi}}, 0\right)$ on $\hat\fC_{\sigma, \rho, T,r}$.
The same conjugacy problem, only without the variable $z$, was dealt
with by Sauzin, \cite{Sa} (who referred to this problem as finding
the characteristic vector field), see also \cite{LMS}, \cite{RT}.
The only difference here is the presence of the quantity $z$.
However, as there is no differentiation in $z$ in the
``unperturbed'' vector field $\x_0$, the quantity $z$ enters the
conjugacy problem as a parameter, and hence the resolution of the
conjugacy is solely based on the invertibility of the operator
${\partial\over
\partial s}+\spr{\omega}{{\partial\over
\partial \varphi}}$ on $\hat{\cal C}_{\sigma, \rho, T}$. Thus, the proof that the equation (\ref{pthje})
can be conjugated to (\ref{lhje}) reproduces the proof of Lemma 4.4
in \cite{RT} verbatim; we skip it, referring the reader to the
latter or in fact any of the three above-mentioned papers.

As the rest of the claims of Theorem \ref{mth} have been shown
earlier in section 1 to follow from the claim (ii), this completes
the proof of theorem. $\Box$

\subsection*{5. Proof of Theorem \ref{sst}}
\addtocounter{section}{1}\setcounter{equation}{0}

The proof follows the lines of the proof of the KAM theorem for
semi-infinite cylinders in \cite{RT}, incorporating the dependencies
in the ``hyperbolic'' variable $z$ and can be in a sense compared to
Graff's proof of the KAM theorem for whiskered tori, see e.g.
\cite{Gr}, \cite{Z2}.

Consider the differentiation operators
\begin{equation}
D_\omega = \spr{\omega}{D_\varphi}, \;\;\; \Dlw\, =\, \lambda_0 D_s
+ D_\omega. \label{dhw}
\end{equation}
The standard KAM theory depends on solvability of linear PDEs with
the operator $D_\omega,$ in \cite{RT} the operator $\Dlw$ was dealt
with.

Consider a perturbation of $H_\omega$ as in (\ref{pt}). The
principal step in proving the structural stability of Hamiltonian
the (\ref{tmp}) is establishing the fact that the Hamiltonian
$H_\omega$ is stable infinitesimally. This is done by solving the
homological equation in the functional linearization of the problem
(i.e. vindicating an ``iterative lemma''). The standard Newton's
iteration follows, see \cite{Z1}, \cite{Z2}. Parameter dependencies
and smallness conditions were worked out for the case $m=0$ in
\cite{RT}; the case $m>0$ makes no difference in this respect.
Indeed, the estimates in the series of propositions in the Appendix,
dealing with inversion of the first order differential operators
involved, are in essence the same as they are in the latter paper,
where the reader is referred for technical detail.

The unknown quantities $(S,\br)$ appearing in (\ref{zed1}) exist,
provided that one can solve the following set of equations (with the
norm of the solution not exceeding the norm of the right-hand side
by a factor, polynomial in the key parameters, such as analyticity
loss):

\begin{equation}
\begin{array}{lll}
\left[ \Dlw + \spr{z}{\Lambda D_z} \right]\hat{S}_0&=& - f -
\spr{\omega}{\hat\xi} + \hat c,
\\ \hfill \\
\left[ \Dlw + \spr{z}{\Lambda D_z} \right]  \hat{\br} &=& \g +
D^2_{\p\p}H_\omega(p,q)|_{\p=0} (d\hat{S}_0 + \hat{\boldsymbol \xi)}
+ B(\varphi,x,z)\hat{\br} - \hat{\boldsymbol \lambda}_0 -
\hat{\boldsymbol\Lambda}^T_0 z.
\end{array}
\label{invert}
\end{equation}
The system (\ref{invert}) arises by direct substitution of
(\ref{zed1}) into (\ref{tmp}) and omitting terms which are
$O_2(|S|+|\br|+|V|)$. As far as the notation is concerned,
$\hat{\boldsymbol\xi}=(\hat{\xi},0,0),\,
\hat{\boldsymbol{\lambda}}_0=(0,\hat{\lambda}_0,0)$ are $n+1+m$
constant column-vectors and
$\hat{\boldsymbol{\Lambda}}^T_0=(0,0,\hat{\Lambda}^T_0)$ is a
constant $(n+1+m)\times m$ matrix.

The role of the parameters $\hat{c},\hat\lambda_0,\hat\Lambda_0$ (in
addition to $\hat\xi$) is to ensure solvability of (\ref{invert})
within the framework of propositions in the Appendix, i.e. to
guarantee that the right hand side is in the complement to the
kernel of the operator $\Dlw + \spr{z}{\Lambda D_z}$ on $\bse$ for
the first equation and $\bseg$ for the second one. Equivalently,
after the canonical transformation $\Psi(\hat{\ar},\hat S),$ the
momentum-linear part of the Hamiltonian $H_\omega$ would acquire a
term
\begin{equation}
\hat H_\omega = \hat{c}+\hat\lambda_0e+\spr{z}{\hat\Lambda_0\bar
z}.\label{add}
\end{equation}
The term $B(x,z)\hat{\br}$ can be described as follows. If
$\hat{\br}=(\hat\beta,\hat b,\hat{\flat})$, describing the
transformation of the $(\varphi,x,z)$ variables respectively, then
$B(x,z)\hat{\br}$ contributes to the equation for the quantity
$\flat$ only, where it results in the term
\begin{equation}
\Lambda^T(s)\hat{\flat} +
\hat{b}D_s\Lambda^T(s)z\label{spellb}\end{equation} in the right
hand side.

In order to solve the first equation in (\ref{invert}), $\hat{c}$ is
to be found, depending on the still unknown $\hat{\xi}$, such that
the right hand side, call it $v_{\hat{S}_0}$, have zero
$\varphi$-mean $\avr{v_{\hat{S}_0}}=0$, cf. (\ref{aver}). Recall
that the mean it is taken by setting $(s,z)=(-\infty,0)$.

No matter what $\hat{\xi}$, such $\hat{c}$ clearly exists, so we can
assume that the right hand side of the first equation has zero mean.
Then $\hat{S}_0$ exists, in any space $\bseprime$, with $\fp'<\fp$,
by Proposition \ref{prthree}. Observe that  $\hat{S}_0$ is
independent of $\hat\xi$.

The second equation in (\ref{invert}) comprises three (systems of)
equations: for the quantities $\hat \beta,\hat b$ and $\hat{\flat}$.
First one considers the equation for $\hat\beta$ and finds $\hat\xi$
such that the right-hand side, call it $v_{\hat\beta},$ has zero
$\varphi$-mean, i.e. $\avr{v_{\hat\beta}}=0$. Note that the last
three terms in the second equation in (\ref{invert}) do not appear
in the equation for $\hat\beta$. Hence by the non-degeneracy
assumption,
\[
\hat\xi\,=\,-\, \avr{D^2_{\iota\iota} H_\omega(\p,\q)_{\p=0}}^{-1}
\tilde v_{\hat\beta},
\]
where $\tilde v_{\hat\beta}$ embraces the first $n$ components of
the $n+1+m$ vector $\g+D^2_{\p\p}H_\omega(\p,\q)|_{\p=0}d{S_0},$
member of the space $\in\bsegprime,$ for any $\fp'<\fp$. This also
determines the constant $\hat{c}$ in (\ref{add}).

Furthermore, the (scalar) $\hat b$-component of the second equation
in (\ref{invert}) is resolved by Proposition \ref{prfour}. The
equation is not soluble without the condition (\ref{snr}). (The term
constant $\hat\lambda_0 e$ is the only thing here to be added to
Hamiltonian $H_\omega$, because under condition (\ref{snr})
constants exhaust the kernel of the operator $\Dlw+\spr{z}{\Lambda
D_z}$ on the space $\bsem$).

Eventually, the equation for the quantity $\hat{\flat}$ is solved.
This equation deserves special attention, so let us write it down
explicitly as follows:
\begin{equation}\label{ss}
[\Dlw+\spr{z}{\Lambda D_z}-\Lambda^T]\hat{\flat} = v_{\hat{\flat}}.
\end{equation}
Let
$$\begin{array}{ccc}
 v_{\hat{\flat}}(\varphi,s,z)&=&v_{\hat{\flat},0}(\varphi,z)+v_{\hat{\flat},1}(\varphi,s,z),
 \\ \hfill \\
\hat{\flat}(\varphi,s,z)&=&\hat{\flat}_0(\varphi,z)+
\hat{\flat}_1(\varphi,s,z),\end{array}
$$
where the quantity $v_{\hat{\flat},1}(\varphi,-\infty,z)=0$, and so
it satisfies the estimate (\ref{espvan}). Therefore, the quantity
$\hat{\flat}_1(\varphi,s,z)\in[\bseoneprime]^m$, corresponding to
the right-hand side $v_{\hat{\flat},1}$ exists, by Proposition
\ref{prfive}.

It remains to determine $\hat{\flat}_0(\varphi,z)$. Let
$$\hat{\flat}_0(\varphi,z) = \hat{\flat}_{0,0}(\varphi) +
\spr{z}{\hat{\flat}_{0,1}(\varphi)}+O_2(z;\varphi),$$ do the same
expansion for the right-hand side $v_{\hat{\flat},0}$. Then
$\hat{\flat}_{0,0}$ is found by Proposition \ref{prtwo}, cf.
(\ref{furc}). As for the term $\hat{\flat}_{0,1}(\varphi)$, it is
easy to see that the quantity $\spr{z}{const.}$ is in the kernel of
the operator in square brackets in (\ref{ss}). Hence the quantity
$\hat\Lambda_0$ is introduced to ensure that the right hand side
$v_{\hat{\flat}}$ do not contain a constant multiple of $z$. This
having been done, for all $z$, the right hand side
$\spr{z}{v_{\hat{\flat},0,1}(\varphi)}$, where
$v_{\hat{\flat},0,1}(\varphi)$ has zero mean, can be resolved by
Proposition \ref{prone}, (i).

Finally, the component $v_{\hat{\flat},0,2}=O_2(z;\varphi)$ in the
right-hand side $v_{\hat{\flat}}$ of equation (\ref{ss}) gets taken
care of as follows. Consider a monomial $z_1^{k_1}\ldots z_m^{k_m}
u_{k}(\varphi)$, with $|k|=k_1+\ldots+k_m\geq2$. Under the action of
the operator $\spr{z}{\Lambda_0 D_z}-\Lambda_0^T,$ taking into
account the fact that $\Lambda_0$ is {\em diagonal}, one gets some
$z_1^{k_1}\ldots z_m^{k_m}\tilde\Lambda_0 u_{k}(\varphi)$, where the
constant matrix $\tilde{\Lambda}_0$ is diagonal and by the condition
(\ref{snr}) is such that the real part of each diagonal entry is
strictly positive, bounded away from zero uniformly in
$(k_1,\ldots,k_m)$ by some $\lambda>0$, which may be set equal to,
say one tenth of the infimum in the right-hand side of (\ref{snr}).
Then the equation gets resolved term by term in the same way as is
(\ref{tr}) in Proposition \ref{prtwo}, the bound for the norm being
uniform for all powers of $z$.

The proof of infinitesimal stability of Hamiltonian $H_\omega$ will
be complete after diagonalizing the constant matrix
$\Lambda_0+\hat\Lambda_0$ by the linear transformation $z\rightarrow
Lz$, where $L$ is a constant near-identity matrix, such that
$L^{-1}(\Lambda_0+\hat\Lambda_0)L$ is diagonal. This is possible as
long as $\hat\Lambda_0$ is small enough.

As we have mentioned earlier, this suffices to prove Theorem
\ref{sst}, as one can now switch on the Newton's iteration
procedure, see \cite{Z1}, \cite{Z2}. For the estimates, which would
result in the forthcoming qualitative version of the theorem, with
the smallness condition and parameter dependencies, see \cite{RT}.
$\Box$

\subsubsection*{Quantitative statement of Theorem \ref{sst}}
We now present a quantitative statement of Theorem \ref{sst}. The
qualitative assumptions naturally look somewhat tighter than as
stated in Theorem \ref{sst}.

\renewcommand{\theassumption}{\arabic{assumption}$'$}
\renewcommand{\thetheorem}{\arabic{theorem}$'$}

\begin{assumption}\label{asth1}
Suppose $\exists\fp=(\sigma,T,\rho,{r})>0,$ as well as
$(\mu,\nu):\,0\leq\mu<\nu\ll 1,$ such that $H_\omega\in\bsek$ and in
the perturbation (\ref{pt}) one has
\begin{equation}
f\in\bse,\;\g\in\bseg, \hspace{3mm}|f|_{\fp}\,\leq
\,\mu,\;|\g|_{\fp}\,\leq \,\mu\nu^{-1}. \label{aspt}
\end{equation}
Regarding the terms in the expression (\ref{tmp}) for $H_\omega$,
suppose
\begin{enumerate}
\item $\omega\in\R^n$ satisfies (\ref{dio});
\item  $\exists\lambda>0$, such that $\forall s\in\Pi_{T,\rho}$ the eigenvalues $\lambda_1(s),\ldots,\lambda_m(s)$ of
$\Lambda(s)=\Lambda_0+\Lambda_1(s)$ satisfy ${\displaystyle
\lambda\leq\min(\Re\lambda_j(s))\leq\max(\Re\lambda_j(s))\leq\lambda_0-\lambda,\,j=1,\ldots
m}$;
\item $\Lambda_0=diag(\lambda_{0,1},\ldots,\lambda_{0,m},)$ and  ${\displaystyle \forall
k\in\Z^m_+,\,\left|\sum_{j=1}^m k_j\lambda_{0,j} -
\lambda_0\right|\geq\lambda;}$
\item  for any constant $m\times m$ matrix $\tilde \Lambda$, with $\|\tilde\Lambda\|<\lambda,$ the matrix
$\Lambda_0+\tilde\Lambda$ is diagonalizable;
\item  $\exists\,R,M>0,$ such that
$\forall\,(\tilde{\p},\q)\in\B^{n+1+m}_\kappa\times{\mathfrak
C}_{\mathfrak{p}},$ $\|\langle D^2_{\iota\iota}
O_2(\tilde{\p};\q)\rangle^{-1}\|\leq R^{-1},$
$\|D^2_{\tilde{\p}\tilde{\p}} O_2(\tilde{\p};\q)\|\leq M$.
\end{enumerate}
\end{assumption}
Let $0<\fp'<\fp$. Further without loss of generality assume that the
quantities
$\delta=\sigma-\sigma',\de=|\fp-\fp'|,\lambda,R,M^{-1},|\omega|^{-1}\leq1$.
Theorem \ref{sst} now vindicates the existence of a canonical
transformation $\Psi$ such that $(H_\omega+V)\circ\Psi=H'_\omega,$
where $H'_\omega\in \bsekprime$ satisfies Assumption \ref{asth1}
with slightly modified parameters $\lambda_0',\Lambda'(s),R',M'$.
The quantitative results and parameter relations, cf. Assumption
\ref{asth1}, can be summarized as follows.

\addtocounter{theorem}{-1}

\begin{theorem}
\label{sstprime} Under Assumption \ref{asth1}, let
$\kappa'=\kappa/2$ and
\begin{equation}
\begin{array}{lllllllllll} &\varsigma&=& \inf(\gamma\delta^{\tau_n},\lambda),&
&\eta&=&R\inf(M^{-1}\varsigma\de,\nu).
\end{array}
\label{eeta}
\end{equation}
There exists a constant $C$, depending only on $\psi$, as well as
the quantities $n,\tau_n,\psi,\fp,\kappa$, such that if
\begin{equation}
\mu\;\;\;\leq \;\;\; C^{-2}\eta^2\;\;\lesssim\;\;
(R/M)^2\de^2[\inf(\varsigma,\nu)]^2, \label{munot}
\end{equation}
the following estimates hold:
\begin{equation}
\begin{array}{cccccccc}
|\cal {S}|_{\fp'}&\leq &C\mu \varsigma^{-1}, & |\hat{\br}
|_{\fp'}&\leq & C\mu (\eta\varsigma)^{-1},\\ \hfill \\
\lambda_0^{-1}|\lambda_0'-\lambda_0|&\leq&C\mu (\eta\lambda)^{-1}, &
\lambda_{0,j}^{-1}|\lambda_{0,j}'-\lambda_{0,j}|&\leq&C\mu
(\eta\lambda)^{-1},
\\ \hfill \\
R^{-1}|R'-R|&\leq&C\mu (\eta\varsigma\de)^{-1}, &
M^{-1}|M'-M|&\leq&C\mu (\eta\varsigma\de)^{-1}.
\end{array}
\label{transf}
\end{equation}
\end{theorem}
The smallness condition (\ref{munot}) is essentially the same as it
was in \cite{RT}. It is determined by the estimates in the series of
propositions in the Appendix only. These estimates coincide with the
estimates in the corresponding propositions in the latter reference,
where the resulting estimate, the analog of (\ref{munot}) is
discussed in detail.

\subsection*{6. Appendix}
\setcounter{proposition}{0}
\renewcommand{\theequation}{A.\arabic{equation}}
\renewcommand{\theproposition}{A.\arabic{proposition}}
\setcounter{equation}{0}

The appendix contains a series of propositions necessary to resolve
the infinitesimal conjugacy problem in the proof of Theorem
\ref{sst} in this paper. The corresponding first order linear PDEs
are solved by the method of characteristics; the proofs bear a close
relation to lemmata in Chapter 5 in Zehnder's work \cite{Z2}, where
the reader is directed for extra detail.

The first result is adopted from \cite{RT}. It is based on the classical result
regarding the operator $D_\omega$, due to R\"ussmann, \cite{Ru}. The frequency
$\omega$ is assumed to satisfy (\ref{dio}), although this assumption is necessary only in the
context of the operator $D_\omega$.
\begin{proposition}\label{prone}
\begin{enumerate}
\item For a function $v\in\bst$ with $\avr{v}=0$, the solution of
the equation $D_\omega u=v$ exists in the space $\bstprime$ for any
$\sigma'<\sigma$. If
$\sigma-\sigma'=\delta,\,\varsigma=\gamma\delta^{\tau_n}$, then
\[
|u|_{\sigma'}\,\lesssim\, \varsigma^{-1}|v|_{\sigma}.
\]
\item
Let $\mathfrak{p}=(\sigma,T,\rho)$ and $v\in\bs$, with $\avr{v}=0$.
The solution of the equation ${\displaystyle \Dlw u=v}$ exists in
$\bsprime$ for any $\mathfrak{p}'=(\sigma',T,\rho)$ with
$0<\sigma'<\sigma$. If
$\sigma-\sigma'=\delta,\,\varsigma=\inf(\gamma\delta^{\tau_n},\lambda_0^{-1})$,
then
\[
|u|_{\fp'}\,\lesssim \, \varsigma^{-1}|v|_{\fp}.
\]
\item
For $v\in\bsm$, there exists a real constant $c,\,|c|\lesssim
|v|_{\fp}$, such that the solution of the equation ${\displaystyle
\Dlw u=v-c}$ exists in $\bsmprime$ and  for the same $\varsigma$ as
in (ii) one has
\[
|u|_{\fp'}\,\lesssim \,\varsigma^{-1}|v|_{\fp}.
\]
\end{enumerate}
\end{proposition}

\begin{proposition}
Let $\mathfrak{p}=(\sigma,T,\rho)$ and $v\in[\bs]^m$. Consider the
equation
\begin{equation}
[\Dlw - \Lambda(s)]u=v,\label{e1}
\end{equation}
where the matrix $\Lambda(s)\in[\bs]^{m^2}$ is such that any
eigenvalue of the constant diagonal matrix $\Lambda_0=
\Lambda(-\infty)=diag(\lambda_1,\ldots,\lambda_m)$ satisfies
$0<c\lambda\leq \Re\lambda_j\leq C\lambda\leq\lambda_0-c\lambda$ for
some $c,C>0$.

The solution of (\ref{e1}) exists in $\bs$ and  with $\varsigma =
\inf( \gamma\delta^{\tau_n},\,\lambda),$
\begin{equation}
|u|_{\fp}\,\lesssim \, \varsigma^{-1}|v|_{\fp}.
\label{bound}
\end{equation}
\label{prtwo}
\end{proposition}
{\bf Proof.} The characteristic flow of $\Dlw$ is
$\phi_t(\varphi,s)=(\varphi+\omega t, s+\lambda_0 t)$, which clearly
maps ${\cal C}_{\fp}$ into itself.

Decompose $v(\varphi,s)=v_0(\varphi)+v_1(\varphi,s)$ (with
$v_0\in[\bst]^m$ and $v_1\in[\bsone]^m$) and
$\Lambda(s)=\Lambda_0+\Lambda_1(s)$ in the sense of (\ref{dcmp}).
Seek the solution $u(\varphi,s)=u_0(\varphi)+u_1(\varphi,s)$,
expecting to find $u_0\in[\bst]^m$ and $u_1\in[\bsone]^m$. Then for
$u_0$ we have
\begin{equation}(D_\omega-\Lambda_0)u_0=v_0\label{tr},\end{equation}
while $u_1$ should satisfy
\begin{equation} D_t u_1 - \Lambda_1(\phi_t(s))u_1  =
v_1(\phi_t(\varphi,s)),\label{fast}\end{equation} where $D_t$ means
differentiation along characteristics. The solution of equation
(\ref{tr}) involves no small divisors and exists as long as the
matrix $\Lambda_0$ is non-singular and diagonalizable. It is assumed
that $\Lambda_0$ is diagonal, so the system of equations (\ref{tr})
gets separated into $m$ equations:
\begin{equation}(D_\omega-\lambda_j)u_{0,j}(\varphi)=v_{0,j}(\varphi)\label{trsep},\;\;j=1,\ldots,m.\end{equation}
This results in an obvious bound $|u_0|_\sigma\lesssim
\lambda^{-1}|v_0|_\sigma$, as each individual equation in
(\ref{trsep}) gets solved as the Fourier series with the
coefficients
\begin{equation} u_{0,j,k} = {v_{0,j,k}\over -\lambda_j+i\spr{k}{\omega}},\;\;k\in\Z^n.
\label{furc}\end{equation} Note that the constants $c,C$ get
absorbed into $\lesssim$ symbols.

For equation (\ref{fast}) let $h(\varphi,s,t,t')$ solve the
homogeneous equation
$$
D_t h(\varphi,s,t,t') =
\Lambda_1(\phi_t(\varphi,s))h(\varphi,s,t,t'),
$$ with $h(\varphi,s,t,t)=1$. As $\Lambda_1(\phi_t(\varphi,s))=\Lambda_1(s+\lambda_0 t)$, one concludes
that $h$ does not depend on $\varphi$ and moreover $h(s,t,t')=\tilde h(s+\lambda_0 t,s+\lambda_0 t').$

Moreover, for $t'<t<0$  one has the growth condition
\begin{equation} |h(s,t,t')|\lesssim e^{C\lambda(t-t')}.\label{hom1}\end{equation}
 Then, as $|v_1(\varphi+\omega
t',s+\lambda_0 t')|\lesssim e^{s+\lambda_0 t'}|v_1|_{\fp}$, by
definition of the space $\bs$, cf. (\ref{espdec}), the integral in
the right hand side of the representation
\begin{equation}
u_1(\varphi,s,t)=\int_{-\infty}^t h(s,t,t') v_1(\varphi+\omega
t',s+\lambda_0 t')dt' \label{intrep}\end{equation} converges
absolutely for all $t\geq0,$ uniformly in $s$, with the bound
$|u_1|_{\fp}\lesssim (\lambda_0-C\lambda)^{-1}\|v_1\|_{\fp},$ and
$u_1$ will be a member of the space $[\bsone]^m$ as is $v_1$. $\Box$

\begin{proposition}\label{prthree}
Let $\mathfrak{p}=(\sigma,T,\rho,{r})$ and $v\in\bse$, with
$\avr{v}=0$, let $\Lambda(s)$ be such that for all $s$, all its
eigenvalues have positive real parts, bounded from zero by
$\lambda>0$. The solution of the equation
\begin{equation}\displaystyle [\Dlw + \spr{z}{\Lambda D_z}]u=v\label{oper}\end{equation} exists in $\bseprime$ for any
$\mathfrak{p}'=(\sigma',T,\rho,{r}),$ with the bound (\ref{bound})
of Proposition \ref{prtwo}.
\end{proposition}
{\bf Proof.} The characteristic flow of the operator in square
brackets in (\ref{oper}) is $\phi_t(\varphi,s,z)=(\varphi+\omega
t,\,s+\lambda_0t,\,\zeta(z,s,t)),$ where $\zeta(z,s,0)=z$ and
$\dot\zeta=\Lambda^T(s+\lambda_0 t)\zeta.$

By positivity of $\lambda_0$ and the assumption on the eigenvalues of
$\Lambda$, bounded in terms of $\lambda$, the characteristic flow $\phi_t$ is well defined on
$(-\infty,0]\times {\mathfrak C}_{\fp}$, and we have an estimate
\begin{equation}
|\zeta(t')|\lesssim e^{-\lambda|t-t'|}|\zeta(t)|,\;\;t'<t\leq0.
\label{etaest}
\end{equation}
After the decomposition $v=v_0(\varphi,s)+ \spr{z}
{v_1(\varphi,s,z)}$ and the same for $u$, the quantity $u_0$ is
found after Proposition \ref{prone}.

Furthermore, $v_1\in[\bse]^m$ still satisfies $|v_1|_{\fp}\lesssim|v|_{fp}$
(recall that by convention ${r}$ is absorbed into $\lesssim$ symbols),
besides
\begin{equation}(\Dlw + \spr{z}{\Lambda D_z})u_1 + \Lambda u_1=v_1.\label{first}\end{equation}
Now let $h(\varphi,s,t,t')$ solve the homogeneous equation
$$
D_t h(\varphi,s,z,t,t') = -\Lambda
(\phi_t(\varphi,s))h(\varphi,s,z,t,t'),
$$
with $h(\varphi,s,z,t,t)=1$, where $D_t$ is differentiation along
characteristics. Clearly for $t'<t<0$  one has
\begin{equation} |h(\varphi,s,z,t,t')|\lesssim e^{\lambda(t'-t)},\label{hom2}\end{equation}
cf. (\ref{hom1}).

Hence one can let
\begin{equation}
u_1(\varphi,s,z)=\int_{-\infty}^0
h(\varphi,s,z,0,t)v_1(\phi_t(\varphi,s,z))dt, \label{ee}
\end{equation}
which guarantees that
$u_1(\varphi,s,z)_{\fp}\lesssim\lambda^{-1}\|v_1\|_{\fp}$ as well as
the fact that $u_1\in[\bse]^m$. $\Box$

The following Proposition follows immediately from Propositions \ref{prtwo} and
\ref{prthree}.

\begin{proposition}\label{prfive}
Let $\mathfrak{p}=(\sigma,T,\rho,{r})$ and $v\in[\bseone]^m$.
Consider the equation
\begin{equation}
[\Dlw +\spr{z}{\Lambda(s)D_z} - \Lambda^T(s)]u=v,\label{e1l}
\end{equation}
where the matrix $\Lambda(s)\in[\bs]^{m^2}$ is such that for all $s$ all its eigenvalues
have positive real parts, bounded from above by $C\lambda\leq\lambda_0-\lambda$, for some
$\lambda>0$.

The solution of (\ref{e1l}) exists in $[\bseone]^m$, with
\begin{equation}
|u|_{\fp}\,\lesssim \, \lambda^{-1}|v|_{\fp}.
\label{boundd}
\end{equation}
\end{proposition}
{\bf Proof:} The characteristic flow $\phi_t(\varphi,s,z)=
(\varphi+\omega t,s+\lambda_0t, \zeta(z,s,t))$, with
$\zeta(z,s,0)=z,$ of the operator $D_t$ clearly maps ${\mathfrak
C}_{\fp}$ into itself. The solution of the homogeneous equation
$h(\varphi,s,z,t,t')$ satisfies estimate (\ref{hom1}), so by
(\ref{espvan}) it becomes possible to define
\begin{equation}
u(\varphi,s,z,t)=\int_{-\infty}^t h(\varphi,s,z,t,t')
v(\varphi+\omega t',s+\lambda_0 t', \zeta(z,s,t'))dt',
\label{intrepp}\end{equation} which satisfies (\ref{bound}). $\Box$

\begin{proposition}\label{prfour}
Let $\mathfrak{p}=(\sigma,T,\rho,{r})$, consider  equation
(\ref{oper}), with $v\in\bsem$. Suppose $\Lambda(s)\in[\bs]^{m^2}$
is such that all its eigenvalues have positive real part, bounded
away from zero by some $\lambda>0$. In addition, suppose
$\Lambda(-\infty)=\Lambda_0=diag(\lambda_1,\ldots,\lambda_m)$ and
for any $k\in\Z_+^m$ one has
\begin{equation}
|\lambda_0 - \sum_{j=1}^m k_j\lambda_j|\geq\lambda.
\label{hypres}
\end{equation}
There exists a real constant $c,\,|c|\lesssim |v|_{\fp}$, such that
the solution of (\ref{oper}) exists in $\bsemprime$ for any
$\mathfrak{p}'=(\sigma',T,\rho,{r})$ with the bound (\ref{bound}).
\end{proposition}
{\bf Proof.} The variation from Proposition \ref{prthree}) (the
characteristic flow is the same, $D_t$ standing for differentiation
along characteristics) is clearly the fact that
\begin{equation}v(\varphi,s,z)={v_0(\varphi,z)\over
\chi(s)}+v_1(\varphi,s,z),\;\;v_1\in\bse.\label{vdc}
\end{equation} So $u$ also has to
have a term ${u_0(\varphi,z)\over \chi(s)}$. Substituting this term into
(\ref{oper}) we get
\[
{1\over\chi(s+t)} \left( -\lambda_0{D\ln\chi(s+t)} + D_t\right) u_0(\phi_t(\varphi,z))=
{1\over\chi(s+t)} v_0(\phi_t(\varphi,z)).
\]
Note that one can represent $D\ln\chi(s)=1+\chi(s)w(s)$, with $w(s)\in\bs$, so the problem will reduce to
Proposition \ref{prthree} if
we can solve the equation
\begin{equation}
(-\lambda_0 + D_\omega + \spr{z}{\Lambda D_z}) u_0(\phi_t(\varphi,z)) = v_0(\phi_t(\varphi,z)).
\label{extra}
\end{equation}
Try $u_0$ as a monomial $u_{0,k_1,\ldots,k_m}(\varphi)
z_1^{k_1}\ldots z_m^{k_m}$, with $k\in\Z^m_+$, substitute it in the
latter equation, with the monomial $v_{0,k_1,\ldots,k_m}(\varphi)
z_1^{k_1}\ldots z_m^{k_m}$ in the right-hand side. This yields

$$
\left(D_\omega+ \sum_{j=1}^m k_j\lambda_j -\lambda_0\right) u_{0,k_1,\ldots,k_m}(\varphi) =
v_{0,k_1,\ldots,k_m}(\varphi),
$$ which implies $|u_0|_{\fp}\lesssim \lambda^{-1}|v_0|_{\fp}$, by (\ref{hypres}),
cf. (\ref{furc}).

The equation for $u_1$ with the right-hand side $v_1$ from (\ref{vdc}) is now
amenable to Proposition \ref{prthree}, the right hand side being
$\tilde{v}_1= v_1+\lambda_0w(s)u_0(\varphi,z)$. In general
$\avr{v_1}\neq 0$ and should be compensated by the constant $c$; it is not difficult to show that in fact
$c=D^2\psi(0)\avr{v_0}+\avr{v_1}$ (see \cite{RT}, Proposition B.4). $\Box$

\vspace{.2in} \medskip \noindent{\bf Acknowledgement:} Research
supported by EPSRC grant GR/S13682/01.
\bibliographystyle{unsrt}

\end{document}